\documentclass[11 pt]{article}
\usepackage{amssymb}

\usepackage{latexsym}

\usepackage{mathptmx}  
\usepackage{mathtools}

\usepackage{lscape} 
\usepackage{amsmath}
\usepackage{a4wide}
\usepackage{color}
\usepackage{epsfig}
\usepackage{float}
\usepackage{tikz}
\usepackage[labelfont=bf]{caption}
\usepackage{subcaption}
\usepackage{graphicx}
\usepackage{array}
\usepackage[font=small,labelfont=bf]{caption}

\usepackage[left=2cm,top=2cm,right=2cm]{geometry}

\newtheorem{theorem}{\bf Theorem}[section]
\newtheorem{definition}[theorem]{\bf Definition}

\newcommand{\beano}{\begin{eqnarray*}}
	\newcommand{\eeano}{\end{eqnarray*}}

\date{}

\begin{document}		
	\title{\bf Study of higher-order interactions in unweighted, undirected networks using persistent homology}
\author{Udit Raj$^{1}$, Slobodan Maleti\'{c}$^{2}$, Sudeepto Bhattacharya$^{1*}$}
\date{%
	\small{$^{1}$ Department of Mathematics, Shiv Nadar Institution of Eminence (Deemed to be University), Delhi-NCR\\
	$^{2}$ Vinča Institute of Nuclear Sciences-National Institute of the Republic of Serbia, University of Belgrade, Serbia
	} \\
*Corresponding author: Sudeepto.bhattacharya@snu.edu.in
\\[2ex]%
}
\maketitle

\begin{center}
	\textbf{Abstract} 
\end{center}	
Persistent homology has been studied to better understand the structural properties and topological features of weighted networks. It can reveal hidden layers of information about the higher-order structures formed by non-pairwise interactions in a network. Studying of higher-order interactions (HoIs) of a system provides a more comprehensive understanding of the complex system; moreover, it is a more precise depiction of the system as many complex systems, such as ecological systems and biological systems, etc., demonstrate HoIs. In this study, the weighted simplicial adjacency matrix has been constructed using the concept of adjacency strength of simplices in a clique complex obtained from an unweighted, undirected network. This weighted simplicial adjacency matrix is thus used to calculate the global measure, which is called generalised weighted betweenness centrality, which further helps us in calculating the persistent homology on the given simplicial complex by constructing a filtration on it. Moreover, a local measure called maximal generalised degree centrality has also been established for a better understanding of the network topology of the studied simplicial complex. All the generalizations given in this work can be reduced to the graph-theoretic case. i.e., for a simplicial complex of dimension 1. Three different filtration schemes for constructing the sequence of simplicial complexes have been given with the help of both global and local measures, and by using these measures, the topology of higher-order structures of the studied network due to the interactions of their vertices has been compared. Further, the illustration of established definitions has been given using a real-life network by calculating Betti numbers up to dimension two.\\

\textbf{Keywords}: Persistent homology, Weighted simplicial adjacency matrix, Simplicial complex, Higher-order structures.

\textbf{MSC}: 05C82, 55U10, 55U05
\pagebreak

\section{Introduction} 
Interactions in the real-world systems are successfully represented by the graph-theoretic networks, where edges signify the presence of interactions between the units of the system. Networks can be used to model a wide range of complex system, including biological systems, social systems, climate systems, and many more \cite{18,19}. However, these systems may also admit higher-order interactions, which are the interaction between more than two vertices, that cannot be captured by using the graph-theoretic approaches \cite{Boccaletti2023,16,17}. Hence, a mathematical object from combinatorial algebraic topology has been introduced called a simplicial complex that can capture the higher-order as well as the lower-order interactions (interactions between two vertices) present in a system. A simplicial complex exhibits simplices of different dimensions corresponding to the interactions present in the system, which makes it a suitable object for representing a complex system. For example, to represent an interaction between three vertices $v_0$, $v_1$ and $v_2$, it contain a simplex represented as $[v_0, v_1, v_2]$, this simplex has dimension $2$ as every simplex which contains $n$ vertices $[v_0, \ldots, v_{n-1}]$ has dimension $n-1$. A simplicial complex can be defined as a collection of simplices with the downward closure property (every subset of a simplex is also a simplex).

For understanding the topological features and structural study of complex systems, homology has been used for a long time. However, the study of homology on a finite set of points, which involves the calculations of homology groups on a simplicial complex, gives us the trivial results, which further gives us the trivial values for all the Betti numbers $(\beta_n)$ except $\beta_0$. Moreover, it can be understood that using the lens of homology, all the simplicial complexes defined on a finite set are identical. Hence, it is better for us to move on to the sequence of a simplicial complex rather than concentrating on the simplicial complex, which is called filtration, and this process is called persistent homology. Persistent homology (PH) can be viewed as a tool to understand the structure of a given weighted network by constructing and calculating its filtration and Betti numbers, respectively. Persistent homology has been widely applied to real-world complex systems, such as for the study of structure in material science, as well as to analyse sensor networks, etc \cite{22,23}.

In this work, we studied about the persistent homology on a given simplicial complex $C$, which is obtained from an unweighted, undirected network using the clique complex method. Persistent homology (PH) is a tool that is used as a feature extractor for a given complex network and studies the network topology. In this work, Betti numbers $\beta_n$ up to dimension $2$ have been calculated, as homology groups for a finite space are trivial for dimensions greater than the space. PH requires a sequence of simplicial complexes at every step for the calculation of homology groups, which, further, are compared by calculating their Betti numbers $\beta_n$. Converting a simplicial complex into a sequence of sub-complexes such that it becomes a nested sequence of sub-complexes called filtration, assuring that the sub-complex at step $m$ should be contained in the sub-complex at step $n$, where $m\leq n$. For a comprehensive and thorough introduction to the computation of persistent homology, one can read the survey articles \cite{5,6}. In filtration at every stage, new edges, triangles are added, or it is better to say new simplices are added to every stage, which leads to the birth and death of cycles at every step. The filtration of a simplicial complex can be obtained in a number of ways, such as VR filtration, Clique complex filtration, etc. In this work, two filtration which are, maximal generalised degree filtration and generalised weighted betweenness filtration have been established. 

Previously, some generalised centralities have been defined for pure simplicial complexes, but they cannot identify the contribution of simplices in a generic simplicial complex (a simplicial complex that contains facets of different dimensions) \cite{Bian,4,14}. Therefore, two generalised centralities established in this work can be given for any generic simplicial complex $\Delta$, hence, every simplex of $\Delta$ can assign a score. A simplicial complex of dimension $K$ contains simplices of dimension $0$ to $K$, and corresponding to every dimension, one can define an adjacency matrix to understand the complete connectivity of the complex. Corresponding to the dimension of every simplex, we can define a hypothetical structural level denoted as $k$, where $ 0 \leq k \leq K$, and at these $k$-levels, the adjacency between the simplices of dimension $k$ can be discussed. Previously, many adjacency rules and their corresponding matrices have been defined at these structural levels, denoted as $k$ \cite{15,7}. This paper introduces the concept of the ``strength of adjacency between simplices" at every structural level. This concept measures which simplex of dimension $l$ makes two $k$-simplices adjacent. For example, in a given simplicial complex at the $k$-level, two $k$-simplices $\tau$ and $\gamma$ can be adjacent because of a $k+1$ simplex, $k+2$ simplex, etc., provided it contains both the $k$-simplices as its face. So, the strength of adjacency, which is represented by $S(\tau, \gamma)$, is given as the dimension of that maximum simplex, which contains both $\tau$ and $\gamma$ as its faces. With the help of the strength of adjacency, the weighted simplicial adjacency matrix $A^{wk}$ at different connectivity levels $k$ has been defined. Each entry $a_{ij}$ of $A^{wk}$ represents the strength of adjacency between simplex $i$ and $j$; moreover, it will be a square matrix of order equal to the number of $k$-simplices in the given complex. These adjacency matrices can fully determine a given simplicial complex; moreover, we can now get a set of weighted simplicial adjacency matrices corresponding to an unweighted network from which the simplicial complex has been obtained. It can also be observed that the concept of strength of adjacency coincides with the standard adjacency rule of simple graphs, since for any two vertices, the strength of adjacency is always one, reflecting the fact that the highest-dimensional simplices that contain them are of dimension one (edges).

After defining the weighted simplicial adjacency matrix, two generalised centrality indices have been defined, which are further used to give weights to the simplices of the studied simplicial complex. The first centrality is maximal generalised degree centrality, which counts the number of maximal simplices (facets) of all the dimensions incident on a given simplex and is given by the summation of their dimension. For example, if two facets of dimensions $2$ and $3$ are incident on vertex $v$, then the maximal generalised centrality of $v$ will be $2+3=5$. This centrality can be calculated for any generic simplicial complex and complements the existing generalised centrality established by Courtney and Bianconi \cite{Bian}. Moreover, the maximal generalised degree centrality can also be reduced to the simplicial complex of dimension one, and can be seen as degree centrality of a graph. Graph-theoretic degree centrality is the number of edges that are incident on a certain vertex. This is the same as adding up the dimensions of all the graph's maximal simplices that contain that vertex as a face, as edges are the maximal simplices of a graph (simplicial complex of dimension one). The second centrality is generalised weighted betweenness centrality, which depends upon the definition of generalised weighted walk and contributes as a global property of the complex. Using these generalised centralities, each simplex will be assigned a weight by its centrality score, and persistent homology can be calculated on them. We argue that by applying the weights, identification of networks with structural similarities can be achieved. Application of persistent homology to a weighted simplicial complex, where weights are associated with generalised
centrality measures, provides insights into the similarity of complex networks that originate from the
hidden higher-order relationships within the structure. Furthermore, it captures the inherent dynamics of multifaceted higher-order structural changes through the filtration process. For studying real world systems at different scales, it is useful to develop a unified mathematical framework
that captures the system’s behavior at different strata. In this sense, and in the context of systems successfully modeled by graphs, the research presented in this paper aims to contribute to establishing generalised centrality measures. The basic criteria for the derivation of these measures are that their limiting cases are standard graph-theoretic counterparts. Hence, the unified mathematical framework for examining complex networks
encompasses both, the simple (pairwise) and higher-order interactions.

The contribution of this article is as follows: (1) It gives the weighted simplicial adjacency matrix $A^{wk}$ of a clique complex at different connectivity levels $k$, using the concept of strength of adjacency. (2) A local measure called maximal generalised degree centrality has been introduced, which can assign scores to each simplex of a given generic simplicial complex. (3) The definition of maximal generalised degree is used to give a local measure for a generic simplicial complex. This measure is called the generalised clustering coefficient.
(4) A global measure called generalised weighted betweenness centrality is also given using the weighted simplicial adjacency matrix. (5) Using the filtration schemes, which depend on global and local measures, persistent homology and Betti numbers up to dimension $2$ are discussed. In the next section, we will discuss some preliminaries required to understand the findings of the paper. 

\section{Preliminaries}
This section provides some standard definitions in combinatorial algebraic topology, to make the article self-contained. Most of the definitions presented in this section can be found in the literature \cite{9,11,12,13,21,Raj}.

\begin{definition}
	Let $V$ be a non-empty, finite set with cardinality $n\in \mathbb{N}$. The elements of $V$ are the vertices, and denoted by $v_i, i=1,\ldots,n$.
	A simplex is a non-empty member of $\mathcal{P}(V)$, here, $\mathcal{P}(V)$ denoted the power set of $V$.
\end{definition}

\begin{definition}
	Given a non-empty set of vertices $V$, of cardinality $n$, family $C$ of simplices is called a simplicial complex if it is closed under subset inclusion, that\,is : \,\,\,
\end{definition}
$$ 
\forall \sigma_i \in C \,\, and \,\, ( \sigma_j\in \mathcal{P}(\mathbb{V}),\,\, if\,\,  \sigma_j\subset\sigma_i  \implies\sigma_j\in C).
$$ \label{11}
By the above condition, the vertex set of the simplicial complex $C$ can be written as $V(C)=\{v\in V |\{v\} \in C\}$. Moreover, taking a subset of simplices of $C$, with the subset inclusion property, is called a sub-complex of the simplicial complex $C$.

\begin{definition}
	Given a simplicial complex $C$, which is a collection of simplices and the dimension of any simplex $\sigma\in C$, $dim(\sigma)$ is defined as $dim(\sigma)=|\sigma|-1$.
\end{definition}

For example, in simplicial complex $C=\{\{1\},\{2\},\{3\},\{1,2\},\{2,3\},
\{3,1\},\{1,2,3\}\}$, the dimension of simplex \{1,2,3\} is $3-1$=$2$. It contains $3$ simplices of dimension $0$, $3$ simplices of dimension $1$, and $1$ simplex of dimension $2$.

\begin{definition}
	The dimension of a simplicial complex $C$, denoted by $dim(C)$, is defined to be largest natural number $r$, wher,e $r\geq 0$, such that $C$ contains a simplex of dimension $r$.
\end{definition}

In the above given example, the dimension of the simplicial complex $C$ is $2$, as the simplex \{1,2,3\} is the highest dimensional simplex of $C$. In the next section, we will see some basic definitions related to homology.

\subsection{Homology of a simplicial complex}
In a graph theoretic network, only the pairwise interactions are captured, and to capture the non-pairwise interactions, we move to the simplicial complex. Whenever a network is given, we typically aim to study its structural and topological properties, and a similar study of associated simplicial complexes can be done using the concept of homology adopted from algebraic topology. Homology measures the cycles in a simplicial complex and computes its corresponding homology groups $H_n$ and Betti numbers, $\beta_n$, which further tells us about the $n$-dimension holes in the simplicial complex.

To understand the homology of a simplicial complex $C$, we first need to understand $k^{th}$-chain groups $C_k(C)$ and some special chains which are called cycles and boundaries. Any $k$-simplex $\{v_0,\ldots,v_k\}$ is an oriented $k$-simplex if it has a fixed orientation, which implies that the ordering of the vertices cannot be changed. On swapping the ordering of two vertices, a minus sign is introduced. For example, if we change the order of two vertices $v_i$ and $v_j$ in $\{v_0,\ldots,v_i,\ldots, v_j,\ldots,v_k\}$ the it become $-\{v_0,\ldots,v_j,\ldots, v_i,\ldots,v_k\}$. Moreover, every given ordering of vertices induces an orientation on the simplex defined on these vertices. If all the simplices of a simplicial complex $C$ have orientation, then it is called an oriented simplicial complex. 

\begin{definition}
	For a given simplicial complex $C$, the $k^{th}$-chain group $C_k(C)$ is the collection of all the $k$-chains in $C$, where a $k$-chain is defined as the sum of a linear combination of all the oriented $k$-simplices of $C$.
\end{definition}
If $\alpha\in C_k(C)$, then $\alpha = a_1\sigma_{1}+\ldots+ a_n\sigma_{n}$. Here, all $a_i's$ are the constants from the field $\mathbb{F}$ and each $\sigma_i's$ are $k$-simplices, for $i\in\{1,2,\ldots,n\}$. One can also see that, the chain group $C_k(C)$ is a vector space over a given field $\mathbb{F}$, with all the $k$-simplices as their basis elements. For example, on putting $k=0$, the chain group $C_0(C)$ contains a linear combination of all the vertices of $C$ as its elements.

A linear map, which is called $k^{th}$-boundary linear map, can be defined on chain groups $C_k(C)$, is given as follows.

\begin{definition}
	A linear map $\delta_k(C)$, is defined for every $k\geq0$ in a given simplicial complex $C$, from the $k^{th}$ chain group $C_k(C)$ to $(k-1)^{th}$-chain group $C_{k-1}(C)$, is given as 
	
	$\delta_k(C): C_k(C) \rightarrow C_{k-1}(C)$
\end{definition}

The above given map sends each basis $k$-chain of $C_k(C)$ to $k-1$-chain in $C_{k-1}(C)$ and it is called $k^{th}$ boundary operator. The one special property of the boundary operator is that composing two consecutive boundary operators gives us a zero map, that is, $\delta_k(C) o \delta_{k+1}(C)$ is a zero map. Which further implies that, for each dimension $k$, the image of $\delta_{k+1}(C)$ lies in  Ker$\delta_k(C)$. In general, every boundary operator $\delta_k$ maps a $k$-simplex to the sum of its alternating $(k-1)$-faces by omitting the $i$th vertex.

$[v_0, v_1,\ldots, v_k]\rightarrow \sum_{i=0}^{k} (-1)^i [v_0,\ldots,v_{i-1}, v_{i+1},\ldots,v_k]$

\begin{definition}
	In a simplicial complex $C$, for every dimension $k\geq 0$, the $k^{th}$ homology group is defined as the $H_k(C)=Ker(\delta_k(C)) / Im(\delta_{k+1}(C))$.
\end{definition}

For every $k\geq0$, with each homology group $H_k(C)$ there associate a $k^{th}$ Betti number denoted by $\beta_k(C)$ is given by the dimension of $k^{th}$ Homology group $H_k(C)$.

On working with a finite set of points, all the Betti numbers $\beta_n$ except for $n=0$ give us trivial results, which can be seen as a limitation of homology. On the other hand, persistent homology, which is an extension of homology, works on the sequence of a simplicial complex starting from a set of vertices. At every stage, new simplices are added depending upon the set threshold $\delta$, This threshold can be chosen in a number of ways; for more information on this one, see the articles \cite{2,1}. The order in which simplices are added in each step determines the topological features of the dataset; hence, choosing a threshold value for filtration will always impact the results.

\begin{definition}
	A filtration of a simplicial complex $C$ defined on a vertex set $V$ is a nested family of simplicial complexes denote by $C^{\delta}_{V}$, where $\delta\in \mathbb{R}$. So that the complex at step $m$ is embedded in the complex at step $n$, $m\leq n$, such that $C^{m}$$\subseteq$ $C^{n}$.
\end{definition}

In the above definition of filtration, $\delta$ can be chosen from a set of values, such as weights of edges, weights of simplices, centrality scores, etc. For the filtration, we have chosen the value of $\delta$ from the scores of two constructed generalised centralities, namely maximal generalised degree centrality and generalised weighted betweenness centrality. In the process of filtration, new holes may appear and some previous holes may disappear, depending on the addition of new simplices. For example, at step $m$, the number of holes was $2$, but at the very next stage, which is $m+1$, the number of holes increased to $4$ in this way, the Betti numbers are changing, which further gives us information about the structure. One can also record the timings of birth and death of the holes in a filtration, with the help of persistence bar codes, which further help to understand the homological features.

\section{Main results}
\textbf{Notations}: Throughout the paper $k$-simplex stands for the simplex of dimension $k$ in a simplicial complex of dimension $K$ and $l$ stands for the $l$=max$\{k+1, k+2, $\ldots$, K\}$, such that $l$-simplex contains the $k$-simplices as their face.

\subsection{Weighted adjacency rule for simplices}

A simplicial complex $\Delta$ of dimension $K$ contains simplices of various dimensions. In graphs, which are considered as $1$-skeleton of simplicial complexes, adjacency has been defined between 0-simplices (vertices) only. On the other hand, in a simplicial complex, adjacency has been defined between various-dimensional simplices at various structure levels $k$. For example, at $k=2$, adjacency will be given for triangles (2-simplices) by certain adjacency rules (Upper adjacency, lower adjacency, or by the combination of both). In this work, we are trying to study the persistent homology on the given unweighted, undirected graph by assigning weights to simplices using a global as well as a local measure. The maximal generalised degree centrality will be the local measure and generalised weighted betweenness centrality will be the global measure. For this, we are required to define a generalised weighted adjacency rule by which at every structure level $k$, we can construct a weighted simplicial adjacency matrix, which will further be used for calculating generalised centralities on them.

In order to achieve the maximal generalised degree and generalised weighted betweenness centrality for simplices of $\Delta$, we need to understand the strength of adjacency between simplices. As two $k$-simplices in $\Delta$ can be adjacent due to the presence of a higher-order simplex $\sigma_l$ which contains both the simplices as its faces, where $l=max\{k+1, k+2, \ldots, K\}$. For example, two $0$-simplices are called adjacent if there exists a $1$-simplex, or any other $l$-simplex which contains both the $0$-simplices as its faces. Hence, depending on the dimension of $l$-simplex, we can the define strength of adjacency.

\begin{definition}
	In a simplicial complex $\Delta$ of dimension $K$, strength of the adjacency between two $k$-simplices $\sigma_i$ and $\sigma_j$ is defined as the dimension of the $l$-simplex $\sigma_m$, which contains $\sigma_i$ and $\sigma_j$ as its faces, where $l=max\{k+1, k+2, \ldots, K\}$. It is denoted by $S(\sigma_i,\sigma_j)$.
\end{definition}  
If there does not exist any $l$-simplex, it implies that the strength of adjacency will be zero. 
On putting $K=1$ in the above definition, which is a graph of dimension $1$, the strength of adjacency for every adjacent pair of vertices will be $1$, as edges are the largest dimensional simplices that contain the vertices as faces.

In the literature, many rules for constructing the adjacency matrix at different connectivity levels $k$ in the simplicial complex $\Delta$ of dimension $K$ have been given. According to the upper adjacency, two $k$-simplices are adjacent to each other if there exists a $(k+1)$-simplex which consist both of them as face. Similarly, two $k$-simplices are lower adjacent to each other if there exists a $(k-1)$-simplex which is the common face of both simplices. Estrada in \cite{15} has defined a general adjacency rule for $k$-simplices, which uses a combination of upper and lower adjacency for construction of adjacency matrix at $k^{th}$-level. In this work, we are giving a new definition of adjacency by using the concept of strength of adjacency between two $k$-simplices. We believe that by this definition, one can distinguish between the adjacency of two $k$-simplices depending on the dimension of the simplex, which allows them to be adjacent. As two $k$-simplices can be upper adjacent provided there exists a larger dimension simplex which consist them as faces, but the dimension of that simplex is not fixed. For any pair of $k$-simplices, a $k+1$-dimensional simplex may exist, or some other higher-dimensional simplex.

\begin{definition}
	In a simplicial complex $\Delta$ of dimension $K$, two $k$-simplices $\sigma_i$ and $\sigma_j$$(\sigma_i\neq  \sigma_j)$ are said to be adjacent if there exists a $l$-simplex $\sigma_m$, which consist both of them as face, where $l=max\{k+1, k+2, \ldots, K\}$, and it is given as:
	
	\[
	[A^{wk}]_{ij} =
	\begin{cases}
		S(\sigma_i, \sigma_j), & \text{if}\quad \sigma_i \text{ is adjacent to } \sigma_j \\
		0, & \text{otherwise}
	\end{cases}
	\]
\end{definition}

In the above definition that gives us the weighted simplicial adjacency at different $k$-levels, on putting $K$=$1$, it gives us the graph-theoretic adjacency between vertices ($0$-simplices). As edges ($1$-simplices) are the largest dimensional simplices of any graph, hence, any entry of matrix $A^{wk}$ can have maximum value $1$ (provided vertices are adjacent) for $k$=$0$.

Some observations of the weighted simplicial adjacency matrix are as follows:

\begin{itemize}
	\item 
	Weighted simplicial adjacency matrix coincide with adjacency matrix of simple graphs if $K$=1.
	\item 
	Every non-zero entry of matrix $A^{wk}$ signifies that corresponding $k$-simplices are adjacent.
	\item 
	At $k^{th}$-level, if $a_{ij}=m$, it implies, $k$-simplices $i$ and $j$ are faces of a $m$-dimensional facet.
\end{itemize}

\subsection{Maximal generalised degree centrality}
Weighted simplicial adjacency matrix quantify the strength of adjacency between simplices. Adjacency of two $k$-dimensional simplices is given by the adjacency rule given in Definition $3.2$. Using this rule, if any simplex of dimension $k$ at $k^{th}$-level is adjacent to another $k$-simplex, then instead of $1$ it has value $m$ in the corresponding entry $[A^{wk}]_{ij}$. Here, $m$ is the dimension of largest simplex which contain both the $k$-simplices as their face. 

Study of persistent homology on an unweighted, undirected network required assigning weights to the network either on edges or on vertices. In this work, we have used a global measure and a local measure to assign weights to simplices of a given simplicial complex. Local measure, which is called ``maximal generalised degree centrality," and global measure, ``generalised weighted betweenness centrality," of weighted simplicial complex have been established in this work. Degree centrality of $0$-simplices (vertices) in a graph theoretic network $G$ is given by the number of incident edges, which are incident on $0$-simplex. Generalisation of this concept for the simplicial complex $\Delta$ of dimension $K$ is to not only count the number of simplices incident on $0$-simplices but also consider the dimension of the incident simplices, moreover, it can be given for any $k$-simplex of $\Delta$. The term maximal in the maximal generalised degree centrality refers to the simplices that are incident on the given simplex; it should always be maximal. Hence, the proposed definition is given below.

\begin{definition} The maximal generalised degree centrality of a simplex $\sigma$ in a simplicial complex $\Delta$ is the sum of the dimensions of all the maximal simplices that have dimension more than $\sigma$ and are incident on $\sigma$. It is denoted by ${D_\Delta}(\sigma)$.		
\end{definition}

It can be written as,
\begin{center}
${D_\Delta}(\sigma)$= $\{\sum dim(\gamma_i) | \sigma\subset \gamma_i \,\, \text{for each} \,\, i \, \text{and each} \,\, \gamma_i \,\, \text{is a maximal simplex of} \,\,\Delta\}$. 
\end{center}

If the simplicial complex has dimension $1$ ($1$-skeleton or graph) then for every $0$-simplex (vertex) the maximal generalised degree centrality ${D_\Delta}(\sigma)$ is same as the degree centrality in graphs. As facets (maximal simplices) which are incident on $0$-simplex will be of dimension $1$ only. That further implies that the maximal generalised degree not only counts the facets which are incident on the $\sigma$, but also takes their dimension into consideration.

For example, in Fig. 1., the maximal generalised degree centrality of vertex $4$ will be, ${D_\Delta}(\{4\})$= 1$\times$1 + 2$\times$1 =3, as, vertex $4$ is incident with two facets of dimension $1$ and $2$ only. On the other hand, generalised degree of simplex \{4,5,6\} is zero, as it is itself a facet and no other facet can incident on it. 

\begin{figure}
	\centering
	\includegraphics[width=.49\linewidth]{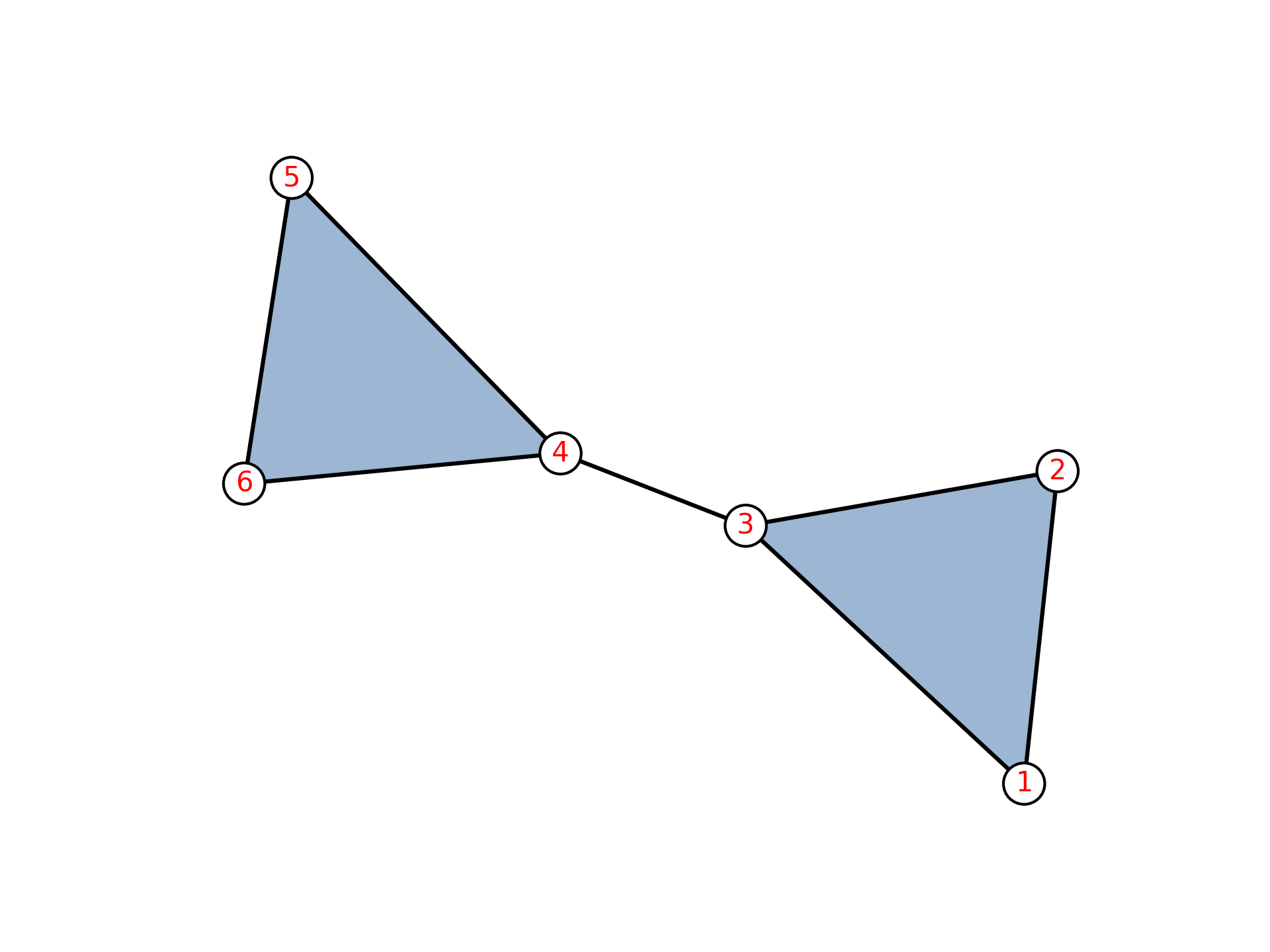}
	\caption{\textbf{Simplicial complex $\Delta$ of dimension 2, it consist facets of dimension 1 and 2 only.}}
	\label{}
\end{figure}
A definition of generalised degree centrality has been introduced by Bianconi in 2015 for the $\delta$ faces of $d$-exclusive simplicial complexes, where $\delta < d$, states that generalised degree centrality of any $\delta$ dimensional face $\sigma_i$ is the number of $d$-simplices incident on $\sigma_i$. But the above definition can be applied to any generic simplicial complex irrespective of it is $d$-exclusive or not. As any generic simplicial complex contains facets of different dimensions and can be incident on a given simplex of $\Delta$. Hence, the definition given in \cite{Bian} can be considered a special case of Definition 3.1. \\

\textbf{Some observations of maximal generalised degree are as follows}:
\begin{itemize}
\item In a simplicial complex $\Delta$ of dimension $K$, maximal generalised degree of any facet is always zero.
\item 
The weighted simplicial adjacency matrix $A^{wk}$ for any $k$ does not behave like the adjacency matrix of simple graphs. For example, the row sum corresponding to any simplex need not give the generalised maximal degree score of that simplex.
\[
\begin{array}{c|cccccc}
	& v_1 & v_2 & v_3 & v_4 & v_5 & v_6 \\
	\hline
	v_1 & 0 & 2 & 2 & 0 & 0 & 0 \\
	v_2 & 2 & 0 & 2 & 0 & 0 & 0 \\
\fbox{$v_3$} & \multicolumn{6}{|c|}{\fbox{$2 \quad 2 \quad 0 \quad 1 \quad 0 \quad 0$}} \\
	v_4 & 0 & 0 & 1 & 0 & 2 & 2 \\
	v_5 & 0 & 0 & 0 & 2 & 0 & 2 \\
	v_6 & 0 & 0 & 0 & 2 & 2 & 0 \\
\end{array}
\]

The above-given adjacency matrix $A^{w0}$ at $k=0$ level is for the simplicial complex $\Delta$, given in Fig.1. It can be seen as an example, where the row sum corresponding to $v_3$ is $5$ and the generalised maximal degree score of $v_3$ is $3$. 
\item 
In a simplicial complex of dimension $1$, $A^{wk}$ coincides with adjacency matrix of graphs.\\

Now, we will see how the generalised clustering coefficient in a simplicial complex can be given using maximal generalised degree centrality.

\subsection{Generalised Clustering Coefficient}

The higher-order clustering coefficient is introduced to thoroughly introduce generalised centrality measures and shape the generalised framework of the analysis of higher-order structures, whose limiting case for $K=1$ is graphs.

The clustering coefficient \cite{Boccaletti2006} is often calculated as a centrality measure that quantifies the density of the neighborhood of nodes in a complex network. Here, we have introduced a measure that transcends the typical clustering coefficient in complex networks and complements the higher-order clustering coefficients defined so far. In this sense, we define the generalised clustering coefficient (GCC) based on the generic definition of the clustering coefficient as a measure of local cohesiveness.

An arbitrary simplex $\sigma_i$ in a given simplicial complex $\Delta$ can have different roles within the structure concerning its immediate neighborhood. 

Due to the richness of the stratified simplicial complex structure \cite{Munkres} and its versatile relations to graphs, different simplicial clustering coefficients have been introduced (see, for example, \cite{Maletic2008}, \cite{Serrano2020}, \cite{Kartun-Giles2019}). Nevertheless, the basic and general motivations for their definition are somewhat similar. 

Compared with the other definitions of the simplicial clustering coefficient, it is clear that they are similar in their intention to relate them to the original complex network clustering coefficient as a special case. Nevertheless, different motivations initiated the derivation of this centrality measure.

The detailed comparison between differently defined simplicial clustering coefficients is outside the scope of the present paper, though comprehensive research on the topic would reveal what kind of information can be deduced from each measure and, furthermore, give insight into how they complement each other.

Generalised clustering coefficient ${C_\Delta}({\sigma_i})$ quantifies the degree of structural cohesion of nearest neighbor simplices of $\sigma_i$.

\begin{definition}
	Generalised clustering coefficient ${C_\Delta}({\sigma_i})$ of a $k$-
	simplex $\sigma_i$ in simplicial complex $\Delta$ is defined by
	\begin{equation}
		{C_\Delta}({\sigma_i})=\frac{\text{$\#$ of adjacences between $\sigma_i$'s
				neighbors}}{\frac{{D_\Delta}({\sigma_i})({D_\Delta}({\sigma_i})-
				1)}{2}}\label{SimplClustCoeff}\, ,
	\end{equation}
	for the adjacency defined in Def. 3.2, and ${D_\Delta}({\sigma_i})$ is
	maximal generalised degree defined in Def. 3.3.
\end{definition}

Instead of generalising the network clustering coefficient, it represents the special case of ${C_\Delta}({\sigma_i})$ when $K=1$, that is simplicial complex is a graph. In other words, if the simplicial complex has dimension 1 (1-skeleton or graph) then for every 0-simplex (i.e. vertex), the generalised clustering coefficient ${C_\Delta}({\sigma})$ is the same as clustering coefficient defined for complex networks.

Hence, in studying the persistent homology, one can assign weights to simplices as a generalised clustering coefficient.

\end{itemize}
\section{Weighted generalised walk and connectivity}
This section establishes a global measure, called generalised weighted betweenness centrality, using the earlier stated weighted adjacency rule in simplicial complexes. Similar to the graph-theoretic weighted networks definition of betweenness centrality, a global measure is required to define weighted generalised walks for weighted simplicial complexes. The definition of generalised walks can be given by the weighted adjacency rule. A weighted generalised walk, in a simplicial complex $\Delta$ of dimension $K$, can be given for the simplices of dimension $k$, the definition is as follows.

\begin{definition}
	A weighted generalised walk between two given $k$-simplices $\sigma_i$ and $\sigma_j$ in a simplicial complex $\Delta$ of dimension $K$ is given by an alternate sequence of $k$-simplices and $l$-simplices such that $\sigma_i=\sigma_1, \sigma_{l_1}, \sigma_2, \ldots, \sigma_{n-1}, \sigma_{l_{n-1}}, \sigma_{n} = \sigma_{j}$ such that for each $i\in\{1, 2, \ldots, {n-1}\}$ $\sigma_i$ and $\sigma_{i+1}$ are adjacent.
\end{definition}
Here, $\sigma_{l_i}$ are the $l$-simplices, where $l=max\{k+1, k+2, \ldots, K\}$, which makes consecutive $\sigma_{i}'s$ adjacent. 
On putting $K$=1 in the above definition, we can get the walk between two $0$-simplices in a simplicial complex of dimension $1$, which is an alternating sequence of $0$-simplices and $1$-simplices or simply the sequence adjacent $0$-simplices, as $\sigma_{l_i}$ which making consecutive $\sigma_{i}'s$ adjacent are of dimension $1$.

The length of a generalised walk $W^{\{wk\}}$ in a simplicial complex $\Delta$ of dimension $K$ between two $k$-simplices is given by the summation of dimensions for all the $l$-simplices in $W^{\{wk\}}$. On putting $K=1$, the length of the generalised walk will be given by the summation of the dimensions of all the $l$-simplices, which is the same as the graph-theoretic case. As the summation of the dimension of $1$-simplices is the same as counting the $1$-simplices.

In Fig.1. the walk between $0$-simplex $6$ and $5$ can be given as, \{6\},\{6,5,4\},\{4\},\{4,3\},\{3\},\{4,3\},\{4\},\{6,5,4\},\{5\}, the length of this walk is $2+1+1+2= 6$, as dimensions of $\sigma_{l_i}$ in the walk is $6$. If we try to minimize the length of this walk then it can be minimized to $2$ as the shortest path between 0-simplices $\{6\}$ and $\{5\}
$ is \{6\},\{6,5,4\},\{5\}.

\begin{definition}
	Shortest path distance between two given $k$-simplices $\sigma_i$ and $\sigma_j$ in a simplicial complex $\Delta$ of dimension $K$ is generalised walk $\sigma_i=\sigma_1, \sigma_{l_1}, \sigma_2, \ldots, \sigma_{n-1}, \sigma_{l_{n-1}}, \sigma_{n} = \sigma_{j}$, in which the length of the generalised walk is minimum.	
\end{definition}

A simplicial complex $\Delta$ of dimension $K$ is said to be connected if there exist a finite length walk between every pair of $0$-simplex in $\Delta$. Only the 0-simplices will decide the connectivity of $\Delta$. Moreover, at every level, one can check the connectivity of the simplicial complex, for example at $k=1$, connectivity can be checked for $1$-simplices and if there exist finite length walk for every possible pair of $1$-simplices, then it is said to be connected at $k$=1 level. Similarly, connectivity can be checked for every level designated by $k$, and $k$ can vary from $0$ to Dim$(\Delta)$. 

\subsection{Generalised weighted betweenness centrality}
In a given simplicial complex $\Delta$ of dimension $K$, the importance of a $k$-simplex can be understood by its capability of lying in the shortest path route of other pairs of $k$-simplices. The generalised weighted betweenness centrality score of a $k$-simplex $\sigma_{i}$ is denoted by $B_{\Delta}(\sigma_{i})$ can be obtained by using the weighted simplicial adjacency matrix $[A^{wk}]_{ij}$ for any $k$. It is called generalised weighted betweenness centrality as it is a generalisation of graph-theoretic centrality and obtained using a weighted simplicial adjacency matrix. 

By using the generalised weighted betweenness score, one can identify the $k$-simplex in a simplicial complex $\Delta$ of dimension $K$, which lies in the shortest path route of maximum $k$-simplices. Moreover, using the scoring of the $k$-simplices filtration of the simplicial complex $\Delta$ can be achieved.

\begin{definition}
	Generalised weighted betweenness centrality $B_{\Delta}(\sigma_i)$ of a $k$-simplex $\sigma_i$ in simplicial complex $\Delta$ is given by  
	\begin{center}
		$B_{\Delta}(\sigma_i)$= $\sum_{\sigma_i \neq \sigma_j \neq \sigma_k}^{}$$\frac{N(\sigma_j, \sigma_k)(\sigma_i)}  {N(\sigma_j, \sigma_k)}$
	\end{center}
\end{definition}
Here, \\
$N(\sigma_j, \sigma_k)$ = Total number of shortest path between k-simplex $\sigma_{j}$ and $\sigma_k$. \\
$N(\sigma_j, \sigma_k)(\sigma_i)$= Total number of shortest path between $\sigma_j$ and $\sigma_k$ which passes through $\sigma_i$.

This is the direct generalization of the graph-theoretic definition, and the above definition can only be applied to $k$-simplices for the calculation of generalised weighted betweenness centrality, provided the given simplicial complex is connected at the $k^{th}$-level. If $\Delta$ is disconnected, then the simplices in the connected component are considered in which $\sigma_i$ lies. Moreover, for the normalization of the scores, we can divide the score by $(n-1)(n-2)/2$, where n is the number of $k$-simplices in the connected component.

Hence, for the calculation of generalised weighted betweenness centrality of $k$-simplex, we first need to find the corresponding weighted simplicial adjacency matrix, then using the algorithm for shortest path we can find $B_{\Delta}(\sigma_{i})$. All the calculations for the generalised weighted betweenness score are done using the networkx package of python, where the shortest path is determined by using the Dijkstra algorithm. 

\subsection{Algorithm for Generalised betweenness centrality}
\begin{itemize}
\item 
Find the weighted simplicial adjacency matrix $A^{wk}$ corresponding to all values of $k$, in which every $a_{ij}$ entry corresponds to weight between simplex $i$ and $j$.
\item 
This weight $w_{ij}$ will be the dimension of the maximum simplex that contains simplex $i$ and $j$ as its face.
\item 
Using the shortest path algorithm of weighted networks, find the number of shortest path for every pair of $k$-simplices. 
\end{itemize}

\subsection{Filtration using the generalised centrality scores}

After defining both the generalised centrality indices, the filtration of the simplicial complex can be achieved by setting up a threshold using the centrality scores. We have defined two types of filtration in this work (A) Maximal generalised degree filtration and (B) Generalised weighted betweenness filtration. In both the filtration, simplices are ranked according to the centrality scores, then a threshold $\delta_m$ is set from the values of centrality scores.\\

\textbf{Maximal generalised filtration}: Let $\Delta$ be a simplicial complex defined on a vertex set $V$, and $D_\Delta$: $\Delta \rightarrow \mathbb{R}$ be a map defined on the simplices of $\Delta$. For any $\delta \in \mathbb{R}$, $\Delta_{\delta}$ is a sub-simplicial complex of $\Delta$ which contains those simplices which have maximal generalised degree score greater than or equal to $\delta$ along with their faces. Such that for any $\delta_m=n_1$ and $\delta_n=n_2$, with $m\leq n$ and $n_2\leq n_1$ then $\Delta_{\delta_m} \subseteq \Delta_{\delta_n}$. The map $D_{\Delta}$ assign a real number, that is maximal generalised degrees score to every simplex of $\Delta$. Hence, after setting a threshold value from the centrality scores, one can get a nested sequence of simplices.\\

\textbf{Generalised weighted betweenness filtration:} Let $\Delta$ be a simplicial complex defined on a vertex set $V$, and $B_\Delta$: $\Delta \rightarrow \mathbb{R}$ be a map defined on the simplices of $\Delta$. For any $\delta \in \mathbb{R}$, $\Delta_\delta$ is a sub-simplicial complex of $\Delta$ which contains those simplices which have generalised weighted betweenness score greater than or equal to $\delta$ along with their faces. Such that for any $\delta_m=n_1$ and $\delta_n=n_2$, with $m\leq n$ and $n_2\leq n_1$ then $\Delta_{\delta_m} \subseteq \Delta_{\delta_n}$.\\

\textbf{Generalised clustering coefficient filtration}: Let $\Delta$ be a simplicial complex defined on a vertex set $V$, and ${C_{\Delta}}: {\Delta}\longrightarrow \Bbb R$ be a map defined on the simplices of $\Delta$. For any $\delta \in \Bbb R$, $\Delta_\delta$ is a sub-simplicial complex of $\Delta$ which contains those simplices which have generalised clustering coefficient score greater than or equal to $\delta$ along with their faces. Such that for any ${\delta_m}={n_1}$ and ${\delta_n}={n_2}$, with $m\leq n$ and ${n_2}\leq {n_1}$ then ${\Delta_{\delta_m}}\subseteq {\Delta_{\delta_n}}$.\\

Using the above-defined filtration schemes, all of them have been used to construct a sequence of simplicial complexes to calculate persistent homology. We have also calculated the various Betti numbers at each sub-complex, allowing us to compare the homology groups at each formed sub-complex. In both filtration schemes, the first sub-complex obtained is not necessarily a set of vertices, but it can contain any simplices that have a centrality score greater than or equal to the maximum score.

\section{Illustrations and Results}
In this section, we illustrate the definitions constructed in this work using a real-world example. The constructed simplicial complex of the Lake network consist $20$ vertices, $51$ $1$-simplices, $24$ $2$-simplices and $5$ $3$-simplices. The data for the same network can be accessed at "An online collection of food-webs, University of Canberra, https://www.globalwebdb.com/, (2012)." The higher-order interactions in this network have been established using the assumed relation of carbon sharing between the units of the network. The clique complex of the lake network has been denoted by $\Delta$, and two filtrations based on the two generalised centrality scores has been constructed for the purpose of Betti numbers comparing. The threshold chosen from maximal generalised degree centrality is denoted by $\delta$, and the threshold chosen from generalised weighted betweenness centrality is denoted by $\delta_B$. Corresponding to every chosen threshold in both cases, a sub-complex has been formed in which simplices that have a score greater than or equal to the threshold value participate. The calculations are done using a python library \textit{Networkx} \cite{3}.

\subsection{Filtration based on maximal generalised centrality scores}

The maximal generalised degree centrality ${D_C}(v_i)$ of vertex $v_{19}$ is the highest as this vertex is incident with $5$ triangles and $5$ tetrahedra. Whereas, $v_{18}$ is incident with only $5$ tetrahedral, which makes its generalised degree score $15$. On comparing the maximal generalised degrees with $1$-skeleton of the Lake network, $v_{19}$ still has highest value. Maximal generalised degree scores for some of the simplices given in Table 5.1.

\begin{table}[H]
	\begin{center}
		\resizebox{9.0cm}{!}
		{\begin{tabular}{|p{1.8cm}|p{4.8cm}|p{3cm}|}
				\hline \textbf{0-simplex $v_i$} & \textbf{Facets of different dimensions incident on $v_i$} & \textbf{Maximal generalised degree centrality ${D_C}(v_i)$}\\
				\hline $v_1$ & 2 simplices of dimension $1$& $1\times2=2$\\
				\hline $v_2$& $6$ simplices of dimension $1$; $3$ simplices of dimension $2$ & $1\times6 + 2\times3= 12$\\
				\hline $v_3$& $4$ simplices of dimension $1$  & $1.4=4$\\
				\hline $v_4$& $3$ simplices of dimension $1$ &$1\times3=3$\\
				\hline $v_5$ &$2$ simplices of dimension $1$ &$1\times2=2$ \\
				\hline $v_{18}$ & $5$ simplices of dimension $3$ & $3\times5=15$ \\
				\hline $v_{19}$ & $5$ simplices of dimension $2$; $5$ simplices of dimension $3$ & $2\times5 + 3\times5= 25$ \\
				\hline
		\end{tabular}} 
	\end{center}
	\caption{{Maximal generalised degree scores ${D_C}(v_i)$ of different 0-simplices in Lake simplicial complex.}}
\end{table}

After calculating the maximal generalised degree scores of every simplex of the Lake network, which ranges from $0$ to $25$, we can now choose a threshold value for the filtration. The first value of threshold is $\delta_{D_1}=25$, which adds all the simplices that have a generalised degree valued greater than or equal to $25$; the corresponding simplicial complex is denoted by $C_{\delta_{D_1}}$. The next threshold which is $\delta_{D_2}=20$ didn't add any new simplex in the previous simplicial complex and the simplicial complex is denoted by $C_{\delta_{D_2}}$. But for the next threshold value $\delta_{D_3}=15$, two new simplices $\{10,18,19\}$ and $\{11\}$ are add up to the sequence, which will change the Betti numbers and other topological properties of the simplicial complex. The Betti number $\beta_0$ for these three threshold values $\delta_1, \delta_2, \delta_3$, changes from $1$ to $2$, which signifies the increase of connected components in the filtration. However, $\beta_1$ and $\beta_2$ remains unchanged, which suggests us the absence of higher-order holes.

For the next threshold values $\delta_{D_4}, \delta_{D_5}, \delta_{D_6}$ and $\delta_{D_7}$, the corresponding simplicial complex adds many simplices, which overall impact the higher-order holes and topology of the simplicial complex. It can be seen in Fig. 5.2. (e) the value of $\beta_1$ is $0$, but for the next simplicial complex in (f), $\beta_1$ increases to 15, which suggests that $1$-dimensional holes are significantly increasing after the threshold value $5$. A few $1$-dimensional holes are disappears in the final simplicial complex given in (g), which give us the value of $\beta_1 = 13$. It suggests that most of the $1$-dimensional holes are persisting after $\delta=2$ and can be seen in the final simplicial complex given Fig. 5.2(g). At $\delta_{D_6} = 2$, the first $2$-dimensional hole appears, whereas at $\delta_{D_7}= 0$, all the $2$-dimensional holes vanish. This indicates that the $2$-simplices responsible for these holes become faces of $3$-simplices (tetrahedra, five in total) that are incorporated into the structure at $\delta_{D_7}= 0$.\\

\textbf{Some interpretation for HoIs}:
\begin{itemize}
	\item 
	In the filtration based on the maximal generalised degree, the first higher-order interaction (HoI) appears at $\delta_{D}= 15$, with additional HoIs emerging at later stages, such as $\delta_{D} = 2$ and $\delta_{D}= 0$. This indicates that the maximal generalised centrality of the HoIs present in the Lake network is nearly zero, with only a few exceptions. Furthermore, this suggests that interactions of order $2$ are not faces of interactions of order $3$ in the Lake network.
	
	\item
	The number of holes of dimension $2$ for threshold value $\delta_{D}=2$ are $5$, which suggest that all the interactions of order $3$ are missing at this stage of filtration and for the last stage of filtration $(\delta_{D}=0)$, $\beta_2$ is zero due the inclusion of all the $3$-simplices.
\end{itemize}

\begin{figure}[H]
	\centering
	\begin{subfigure}[b]{0.2\textwidth}
		\centering
		\includegraphics[width=\textwidth]{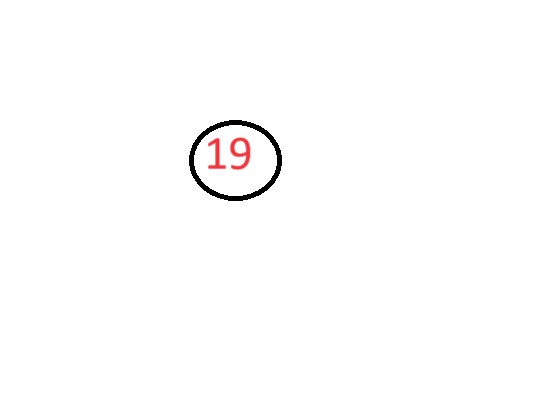}
		\caption{$\delta_{D_1}=25, \beta_0=1$}
	\end{subfigure}
	\begin{subfigure}[b]{0.2\textwidth}
		\centering
		\includegraphics[width=\textwidth]{threshold_25.png}
		\caption{$\delta_{D_2}=20,\beta_0=1$}
	\end{subfigure}
	\begin{subfigure}[b]{0.25\textwidth}
		\centering
		\includegraphics[width=\textwidth]{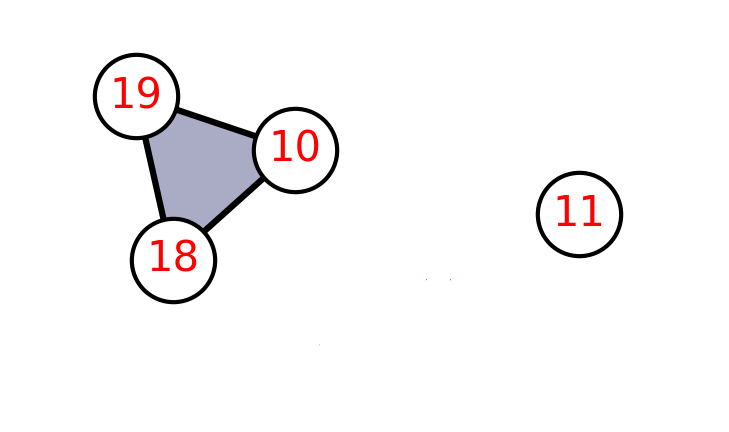}
		\caption{$\delta_{D_3}=15,\beta_0=2$}
	\end{subfigure} \\
	\begin{subfigure}[b]{0.3\textwidth}
		\centering
		\includegraphics[width=\textwidth]{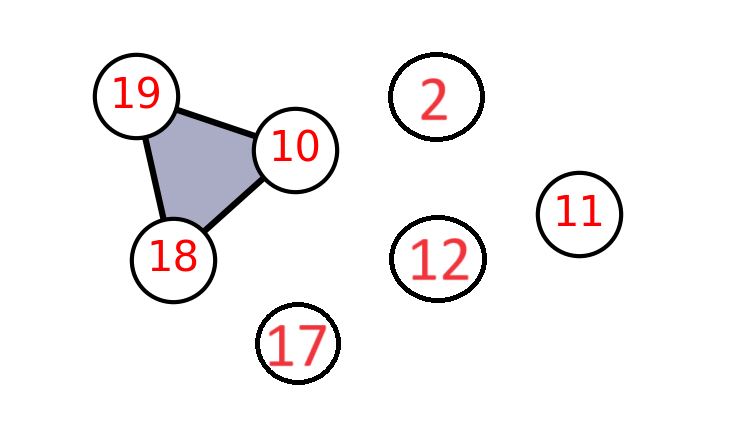}
		\caption{$\delta_{D_4}=10,\beta_0=5$}
	\end{subfigure}
	\begin{subfigure}[b]{0.3\textwidth}
		\centering
		\includegraphics[width=\textwidth]{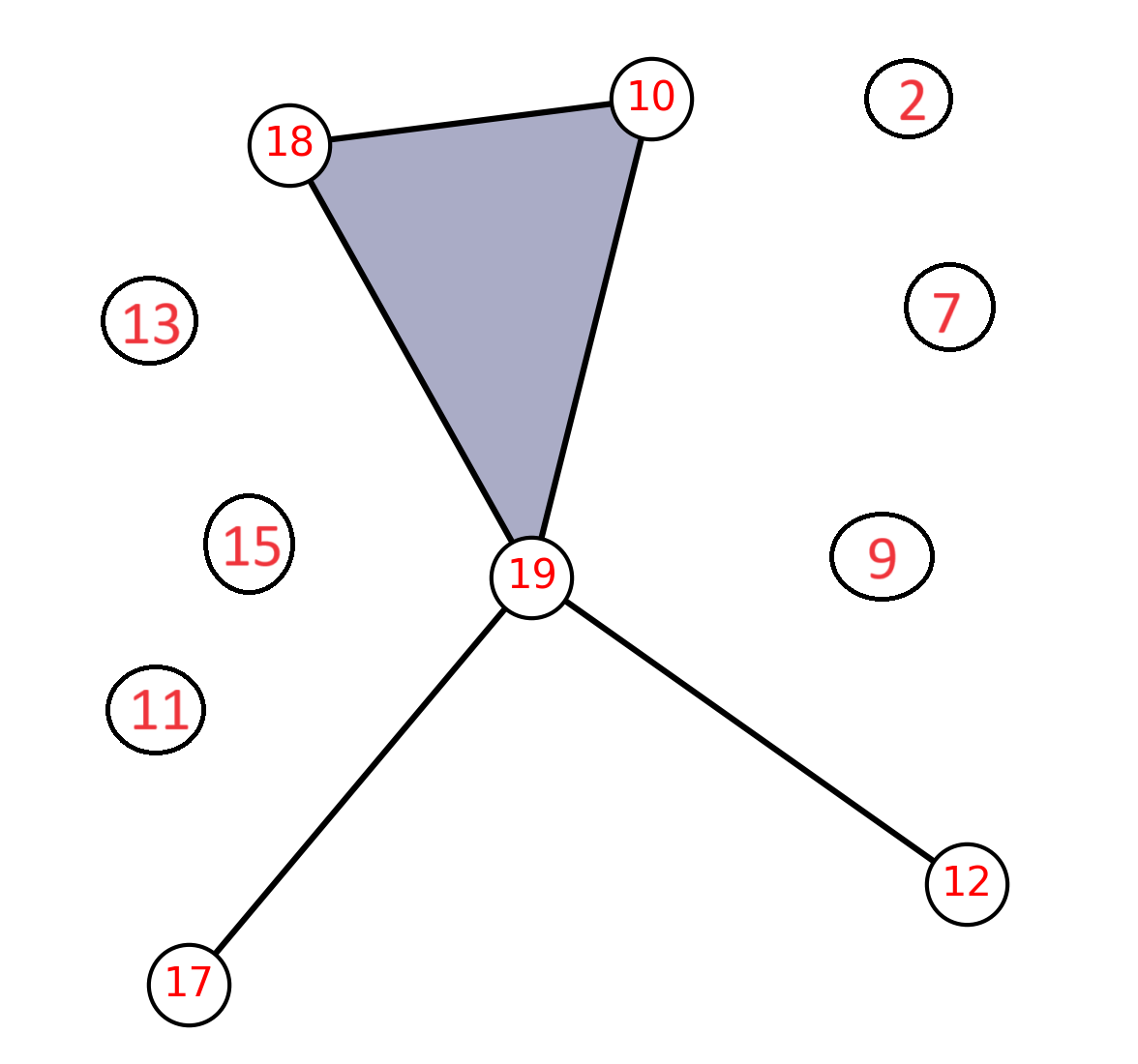}
		\caption{$\delta_{D_5}=5,\beta_0=7$}
	\end{subfigure}
	\begin{subfigure}[b]{0.3\textwidth}
		\centering
		\includegraphics[width=\textwidth]{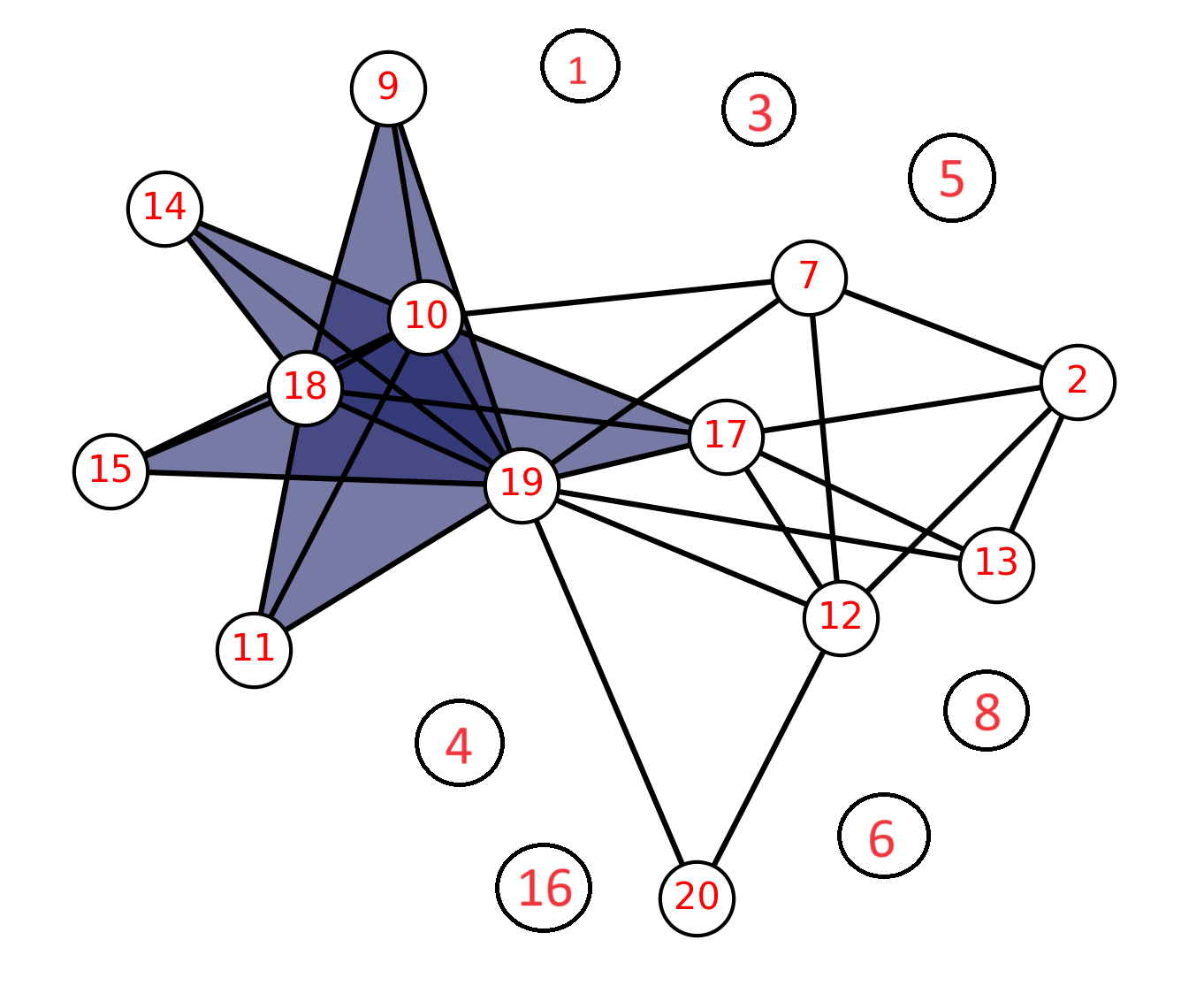}
		\caption{$\delta_{D_6}=2,\beta_0=8$}
	\end{subfigure}
\begin{subfigure}[b]{0.3\textwidth}
	\centering
	\includegraphics[width=\textwidth]{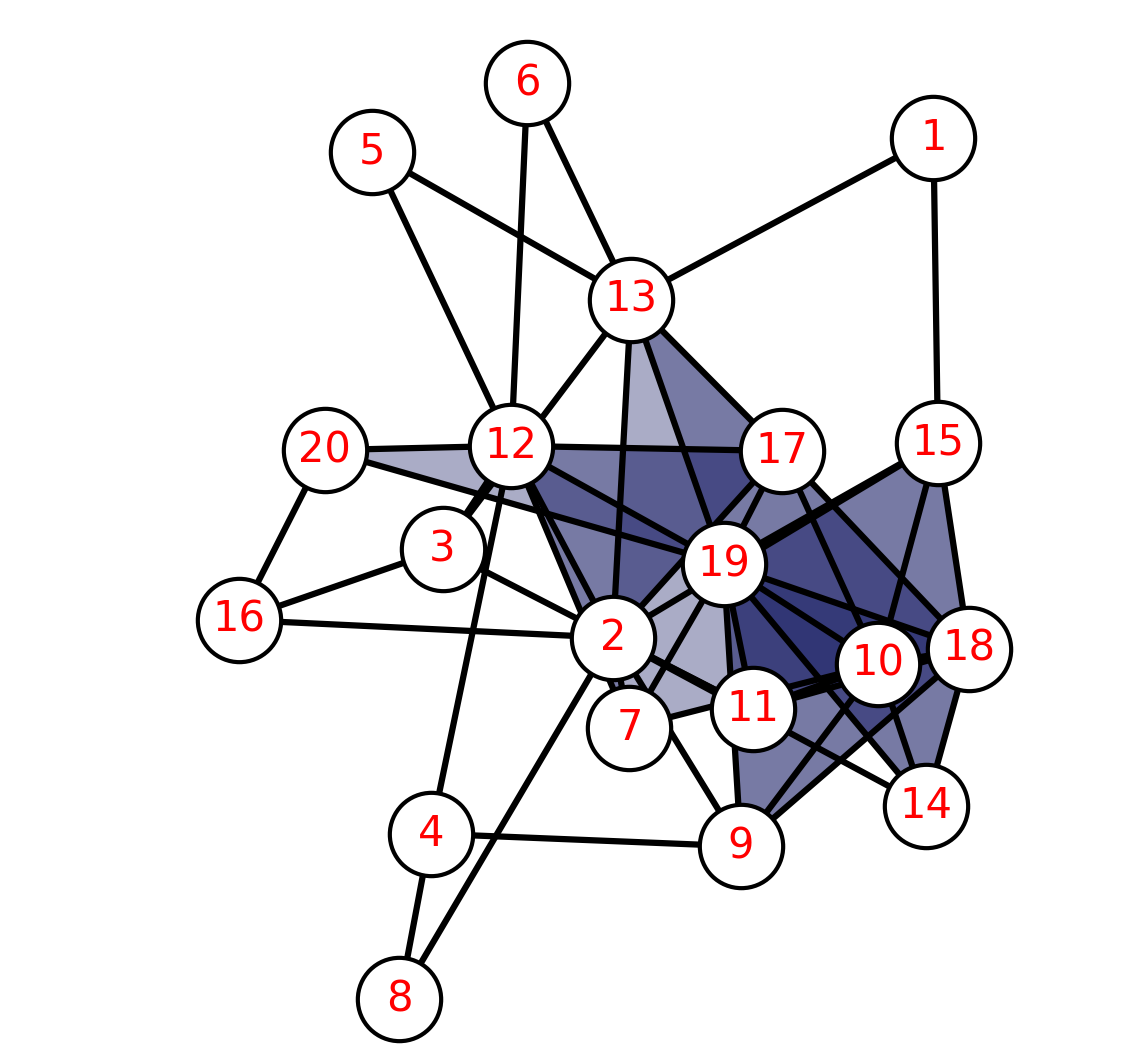}
	\caption{$\delta=0, \beta_0=1$}
\end{subfigure}
	\caption{Filtration of Lake network using the maximal generalised degree scores as threshold for $\delta_{D}={25, 20, 15, 10, 5, 2, 0}$}. 
	\label{fig:six_images}
\end{figure}

\begin{table}[H]
	\begin{center}
		\resizebox{9.1cm}{!}
		{\begin{tabular}{|p{3.8cm}|p{4.8cm}|}
				\hline \textbf{Threshold value for the filtration $\delta_{D_1}$} & \textbf{Betti numbers $\beta_n$} \\
				\hline $\delta_{D_1}=25$ & $[\beta_0, \beta_1, \beta_2]=[1,0,0]$\\
				\hline $\delta_{D_2}=20$ & $[\beta_0, \beta_1, \beta_2]=[1,0,0]$\\
				\hline $\delta_{D_3}=15$ & $[\beta_0, \beta_1, \beta_2]=[2,0,0]$\\			
				\hline $\delta_{D_4}=10$ & $[\beta_0, \beta_1, \beta_2]=[5,0,0]$\\
				\hline $\delta_{D_5}=5$ & $[\beta_0, \beta_1, \beta_2]=[7,0,0]$\\
				\hline $\delta_{D_6}=2$ & $[\beta_0, \beta_1, \beta_2]=[8,15,5]$\\
				\hline $\delta_{D_7}=0$ & $[\beta_0, \beta_1, \beta_2]=[1,13,0]$\\
				\hline
		\end{tabular}} 
	\end{center}
	\caption{{Various Betti numbers corresponding to different threshold values $\delta_D$, based on the maximal generalised degree centrality}}
\end{table}

Calculation of Betti numbers $\beta_n$ for different simplicial complexes formed during the filtration based on maximal degree centrality, give us the understanding of $n$-dimensional holes and other structure details which persist during the filtration.

\subsection{Filtration based on generalised betweenness centrality scores}
On the other hand, the filtration formed by setting the threshold using the generalised weighted betweenness centrality scores shows different values of the initial three Betti numbers, which are given in Table 3. The threshold depending on generalised weighted betweenness centrality is denoted by $\delta_B$ and the range of $\delta_B$ is varies from $0$ to $0.35$ as $0.35$ is the highest value of $B_\delta(\beta)$ for the simplex $\beta=[10,18,19]$ in the lake network. In Table 3, Betti numbers corresponding to different threshold values $\delta_B$ have been given, which suggest that $\beta_0$ is significantly increasing from $\delta_B=0.35$ to $\delta_B=0.02$, and it can be concluded that number of connected components are increasing as we are decreasing the threshold value from $0.35$ to $0$. For $\delta_B=0$, due to the participation of all the simplices $\beta_0$ is 1, which is the original simplicial complex $\Delta$. 
\begin{table}[H]
	\begin{center}
		\resizebox{9.9cm}{!}
		{\begin{tabular}{|p{3.8cm}|p{4.8cm}|}
				\hline \textbf{Threshold value for the filtration $\delta_B$} & \textbf{Betti numbers $\beta_n$} \\
				\hline $\delta_{B_1}=0.35$ & $[\beta_0, \beta_1, \beta_2]=[1,0,0]$\\
				\hline $\delta_{B_2}=0.30$ & $[\beta_0, \beta_1, \beta_2]=[2,0,0]$\\
				\hline $\delta_{B_3}=0.15$ & $[\beta_0, \beta_1, \beta_2]=[3,0,0]$\\			
				\hline $\delta_{B_4}=0.05$ & $[\beta_0, \beta_1, \beta_2]=[7,0,0]$\\
				\hline $\delta_{B_5}=0.02$ & $[\beta_0, \beta_1, \beta_2]=[9,1,0]$\\
				\hline $\delta_{B_6}=0$ & $[\beta_0, \beta_1, \beta_2]=[1,13,5]$\\
				\hline
		\end{tabular}} 
	\end{center}
	\caption{\textbf{Various Betti numbers corresponding to different threshold values $\delta_B$, based on the generalised betweenness centrality scores}}
\end{table}

\begin{figure}[H]
\centering
\begin{subfigure}[b]{0.3\textwidth}
	\centering
	\includegraphics[width=\textwidth]{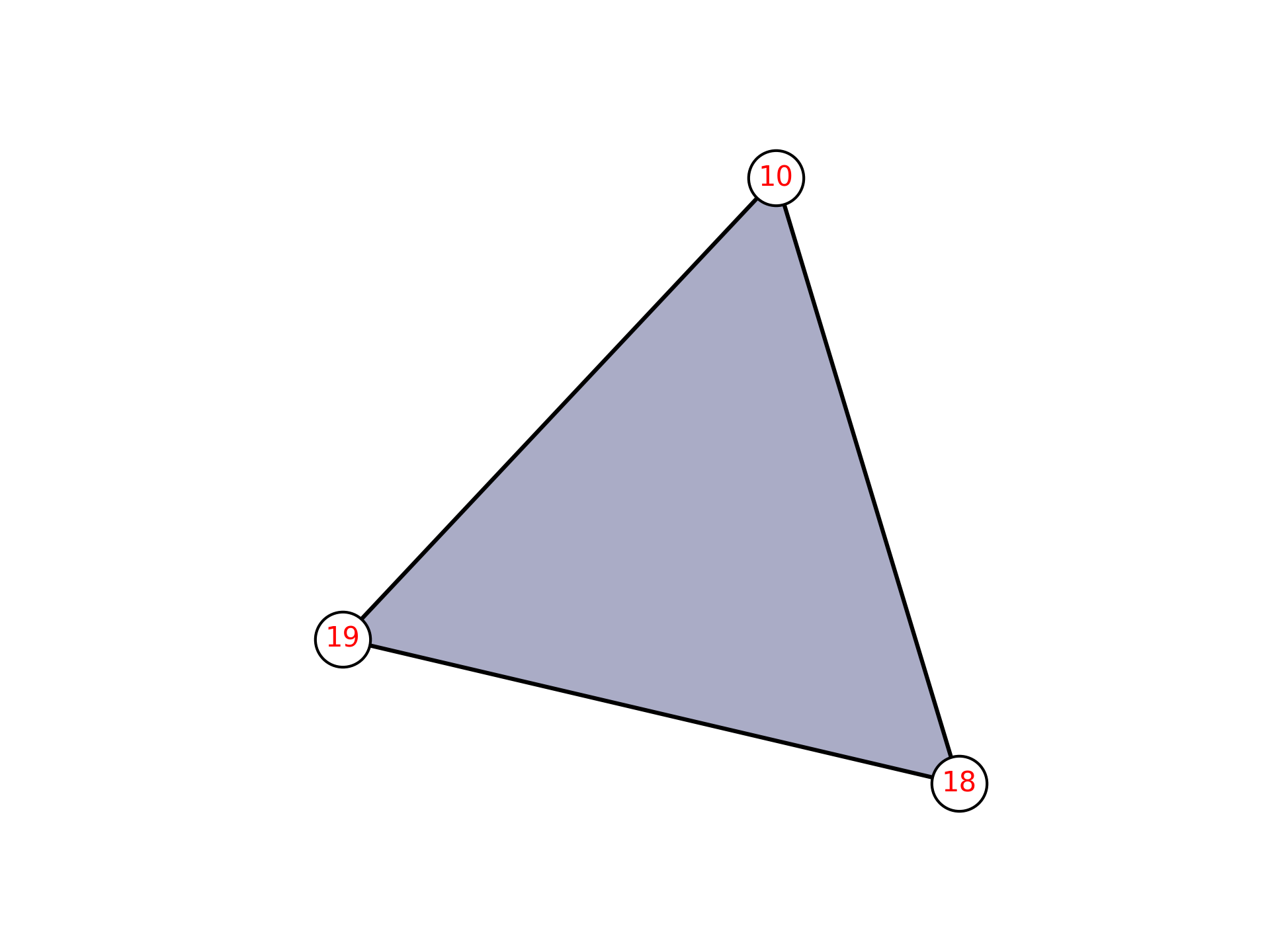}
	\caption{$\delta_B=0.35, \beta_0=1$}
\end{subfigure}
\begin{subfigure}[b]{0.3\textwidth}
	\centering
	\includegraphics[width=\textwidth]{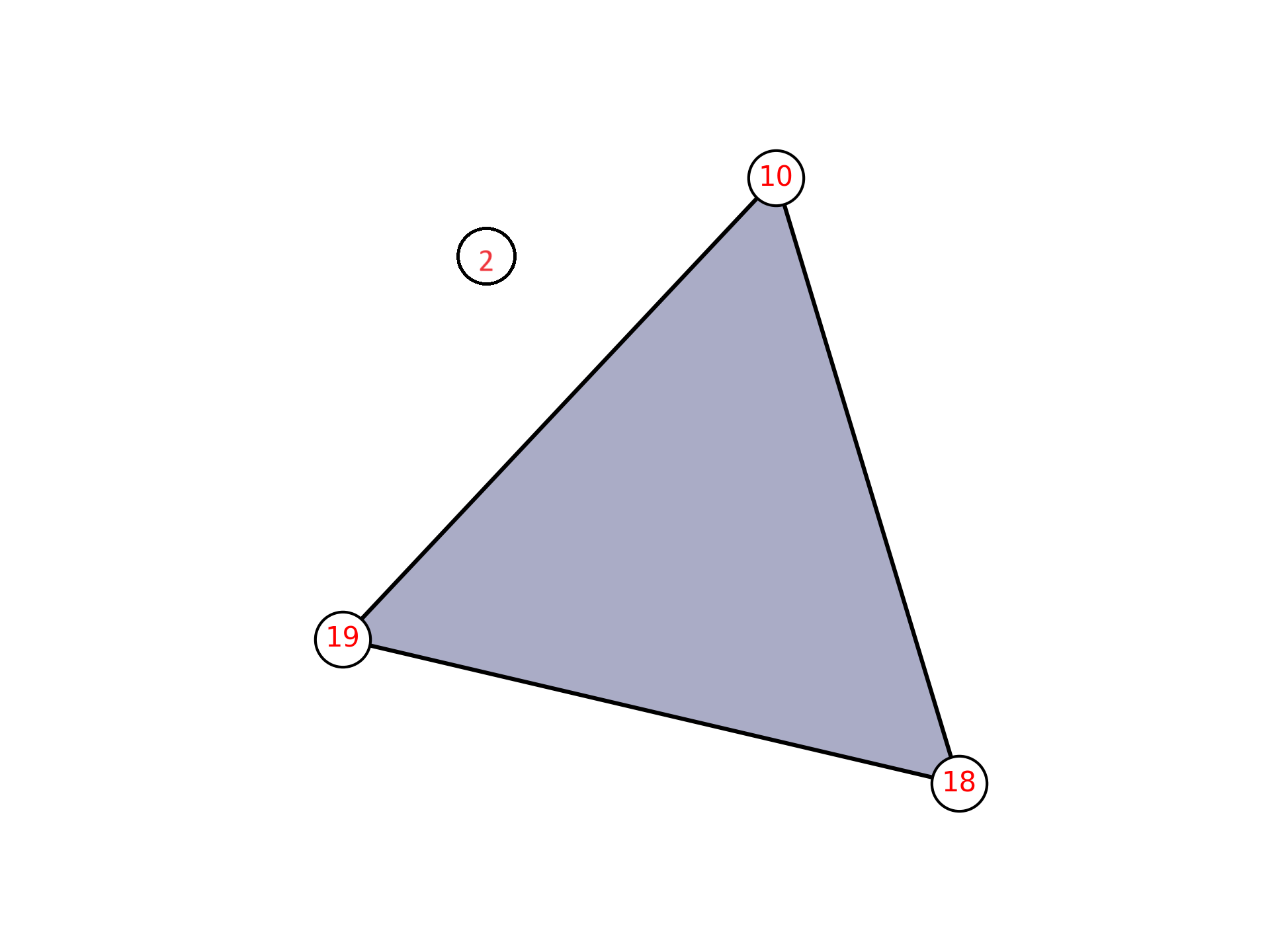}
	\caption{$\delta_B=0.30,\beta_0=2$}
\end{subfigure}
\begin{subfigure}[b]{0.3\textwidth}
	\centering
	\includegraphics[width=\textwidth]{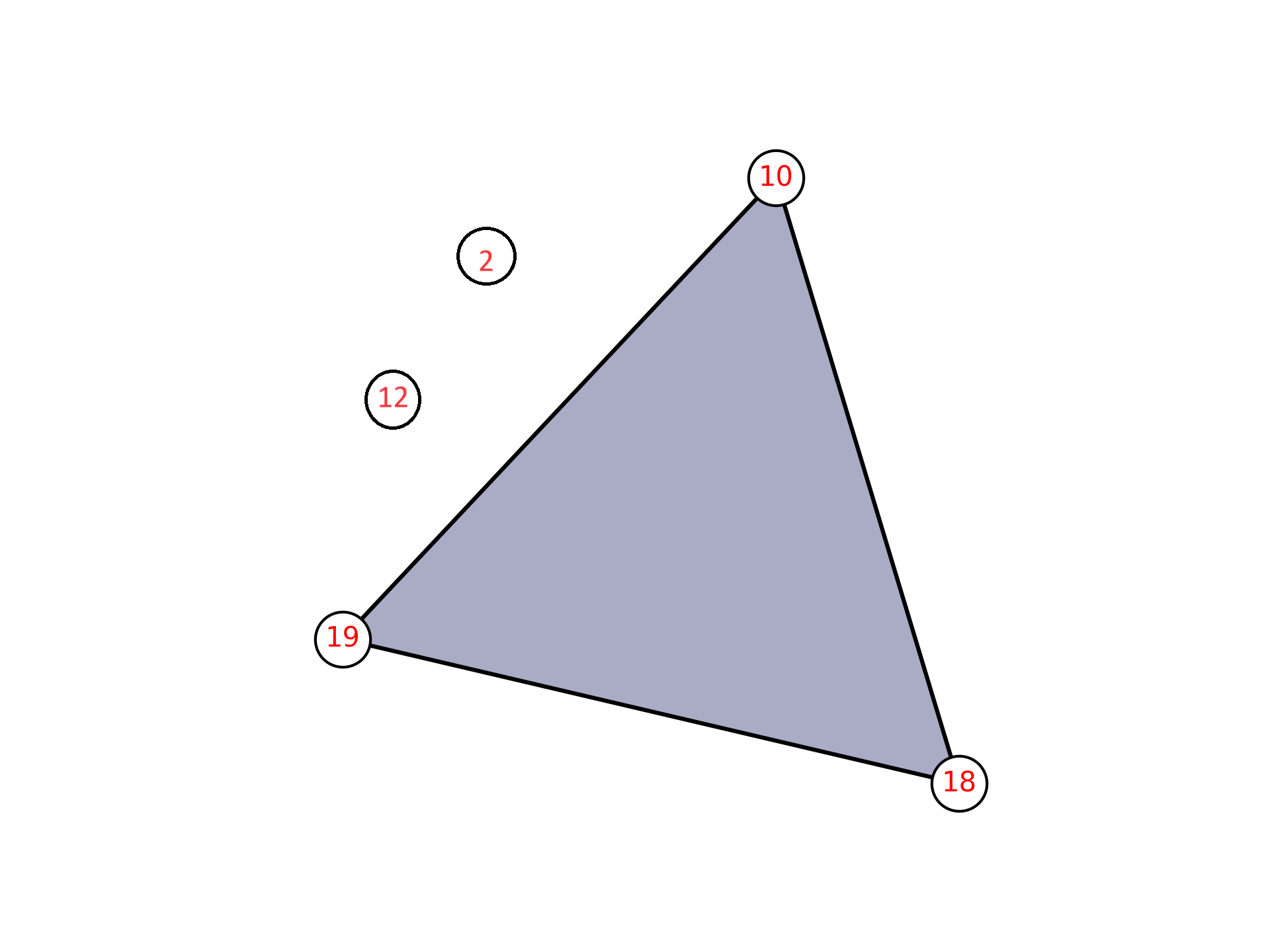}
	\caption{$\delta_B=0.15,\beta_0=3$}
\end{subfigure} \\
\begin{subfigure}[b]{0.3\textwidth}
	\centering
	\includegraphics[width=\textwidth]{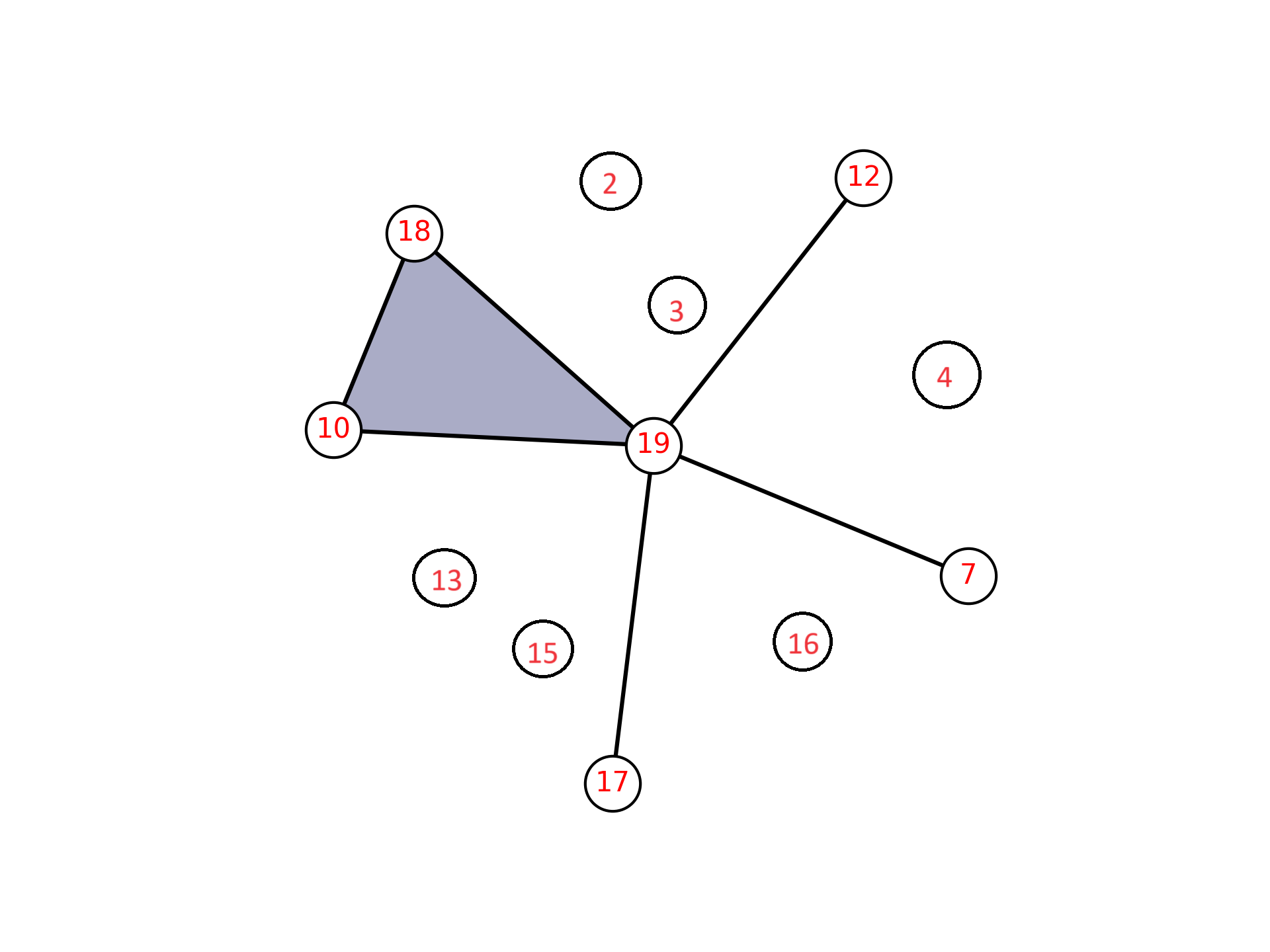}
	\caption{$\delta_B=0.05,\beta_0=7$}
\end{subfigure}
\begin{subfigure}[b]{0.3\textwidth}
	\centering
	\includegraphics[width=\textwidth]{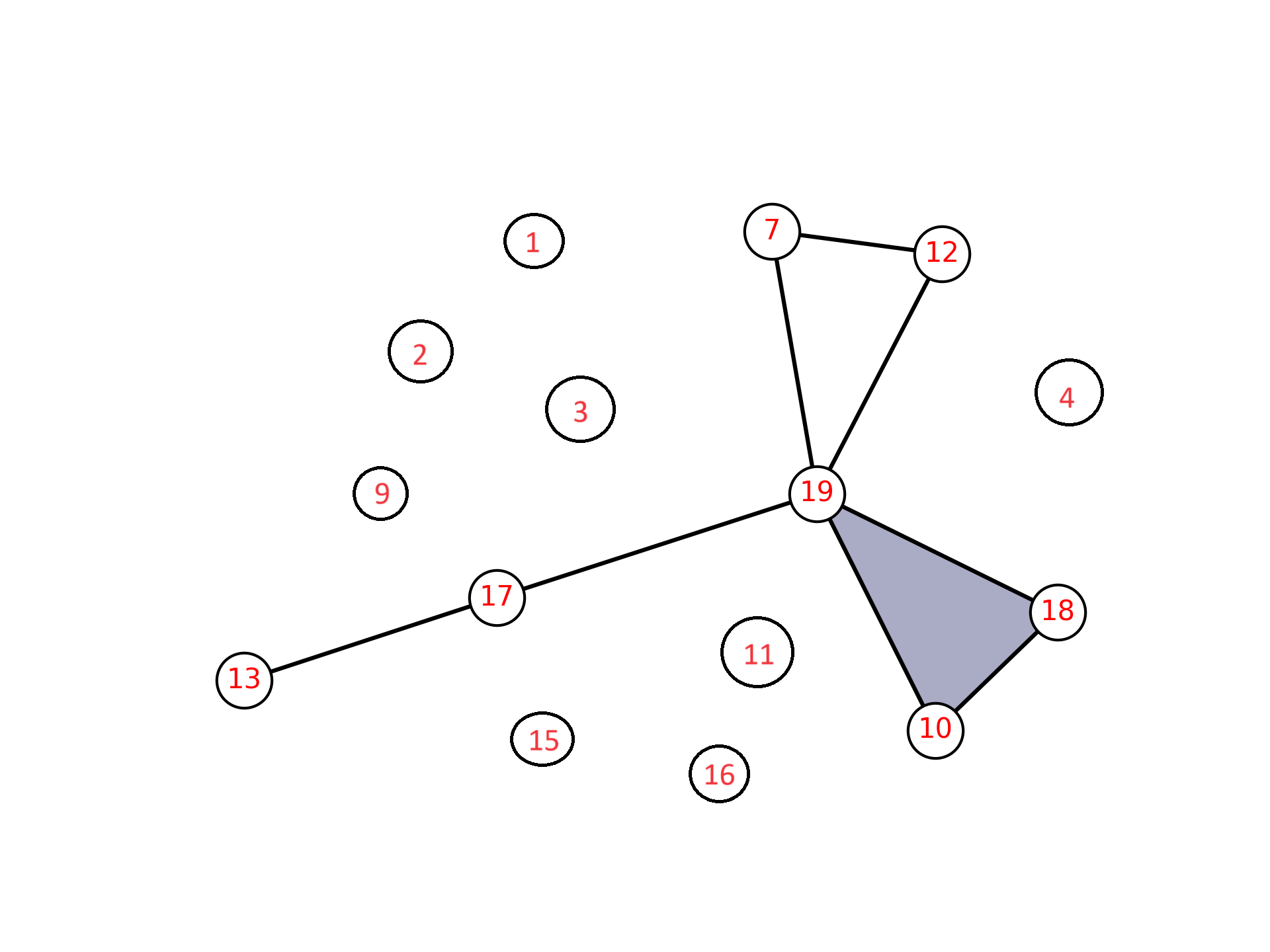}
	\caption{$\delta_B=0.02,\beta_0=9$}
\end{subfigure}
\begin{subfigure}[b]{0.3\textwidth}
	\centering
	\includegraphics[width=\textwidth]{threshold_0.png}
	\caption{$\delta_B=0, \beta_0=1$}
\end{subfigure}
\caption{Filtration of Lake network using the generalised betweenness centrality scores as threshold for $\delta_B={0.35, 0.30, 0.15, 0.05, 0.02, 0}$}. 
\label{fig:six_images}
\end{figure}

The values of $\beta_1$ for all the threshold values is $0$, except for $\delta_B=0.02$ and $\delta_B=0$, this is because the cycles of dimension $1$ are absent in the all the sub-complexes except for the last two. The values of $\beta_1=1$ and $13$ for the last two sub-complexes ('e' and 'f') of Fig.3 suggest that the cardinality of $1$-dimensional holes is $1$ and $13$ respectively in the corresponding sub-complexes. For example, cycle formed by $1$-simplices $[7,19], [7,12]$ and $[19,12]$ is the $1$-dimensional hole in the sub-complex formed for $\delta_B=0.02$. This significant increase in the $1$-dimensional holes suggests that most of the $1$-simplices in the $\Delta$ has generalised weighted betweenness score is between $0.02$ and $0$.\\

\textbf{Some interpretations about HoIs}
\begin{itemize}
\item
No, $2$-dimensional holes appeared in the filtration based on weighted generalised betweenness centrality except for $\delta_B=0$, as all the $2$-simplices that are participating are not forming any solid tetrahedral, which further suggest that these $2$-dimensional cycles are not boundaries of any $3$-simplex. 
\item 
The higher-order interaction $[10,18,19]$, formed at $\delta_B = 0.35$, persists until the end and can be regarded as a dominating higher-order interaction (HoI) in the Lake network.
\item 
The loop that is created at $\delta_B = 0.02$ disappears at the next threshold, $\delta_B = 0$, indicating that this $1$-dimensional hole does not contribute to the final simplicial complex.
\end{itemize}

\subsection{Filtration based on generalised clustering coefficient scores}

In the filtration formed by the generalised clustering coefficient scores of the studied network, we have noticed that the GCC scores for some simplices are more than 1. For example, $1$-simplex $\sigma_{i}=[11,18]$ has GCC value $3.33$, which implies that the denominator that is given by  $\frac{D_\Delta(\sigma_{i})(D_\Delta(\sigma_{i})-1)}{2}$, can have a value less than the number of adjacencies between neighbors of $\sigma_{i}$, as maximal generalised degree counts the number of maximal faces incident on a given simplex; hence, maximal degree cannot decide the maximum number of neighbor adjacencies for a given simplex having dimension greater than 1. To make the value of GCC for each simplex between $0$ and $1$, these values are normalized by dividing the maximum value of GCC for that dimension. For instance, the maximum GCC value for $1$-simplices was $3.33$, leading to the normalization of all values by dividing by 3.33. In the next table, all the Betti numbers ($\beta_0, \beta_1, \beta_2$) corresponding to the filtration formed by setting threshold values $\delta_{CC}$ are given below.

\begin{table}[H]
	\begin{center}
		\resizebox{9.9cm}{!}
		{\begin{tabular}{|p{3.8cm}|p{4.8cm}|}
				\hline \textbf{Threshold value for the filtration $\delta_{CC}$} & \textbf{Betti numbers $\beta_n$} \\
				\hline $\delta_{{CC}_1}=1$ & $[\beta_0, \beta_1, \beta_2]=[1,14,5]$\\
				\hline $\delta_{{CC}_2}=0.66$ & $[\beta_0, \beta_1, \beta_2]=[2,14,5]$\\
				\hline $\delta_{{CC}_3}=0.30$ & $[\beta_0, \beta_1, \beta_2]=[1,6,5]$\\			
				\hline $\delta_{{CC}_4}=0.16$ & $[\beta_0, \beta_1, \beta_2]=[1,7,5]$\\
				\hline $\delta_{{CC}_5}=0.09$ & $[\beta_0, \beta_1, \beta_2]=[1,15,6]$\\
				\hline $\delta_{{CC}_6}=0$ & $[\beta_0, \beta_1, \beta_2]=[1,13,5]$\\
				\hline
		\end{tabular}} 
	\end{center}
	\caption{\textbf{Various Betti numbers corresponding to different threshold values $\delta_{CC}$, based on the generalised clustering coefficient scores.}}
\end{table}

The value of $\beta_0$ for all the threshold values is $1$, except for threshold $0.66$. Hence, throughout the filtration, the number of connected components is $1$ except for the threshold of $0.66$. It is also evident from Table 4 that $\beta_2$ is always $5$ except for the threshold value 0.09; this behavior can be explained by the addition of the 2-simplex [10,18,19], which at threshold 0.09 increases the number of 2-dimensional holes by one.

Six filtration stages are graphically presented in Fig. 4 (a) – (f) for threshold values $\delta_{CC}= \{ 1, 0.66, 0.30, 0.16, 0.09, 0 \}$, respectively.
The filtration due to the clustering of simplices produces the dynamics of emergence of invariant topological objects at different dimensions. Furthermore, the irregularity of Betti number values when the filtration threshold is decreasing suggests that simplices may form nontrivial intrinsic local relationships embedded into the structure imposed by the higher-order clustering. 

\begin{figure}[H]
	\centering
	\begin{subfigure}[b]{0.3\textwidth}
		\centering
		\includegraphics[width=\textwidth]{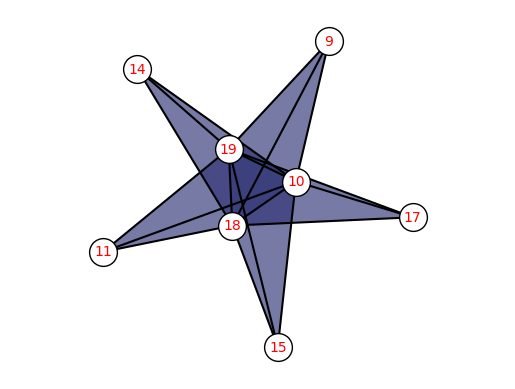}
		\caption{$\delta_{CC}=1, \beta_0=1$}
	\end{subfigure}
	\begin{subfigure}[b]{0.3\textwidth}
		\centering
		\includegraphics[width=\textwidth]{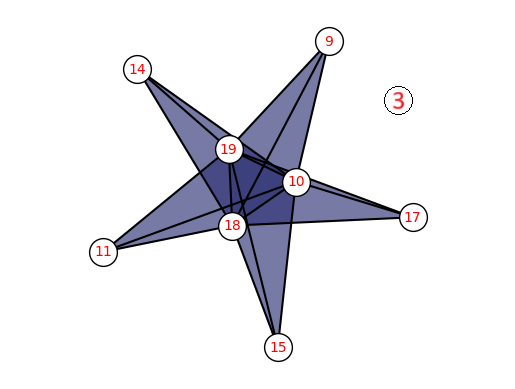}
		\caption{$\delta_{CC}=0.66,\beta_0=2$}
	\end{subfigure}
	\begin{subfigure}[b]{0.3\textwidth}
		\centering
		\includegraphics[width=\textwidth]{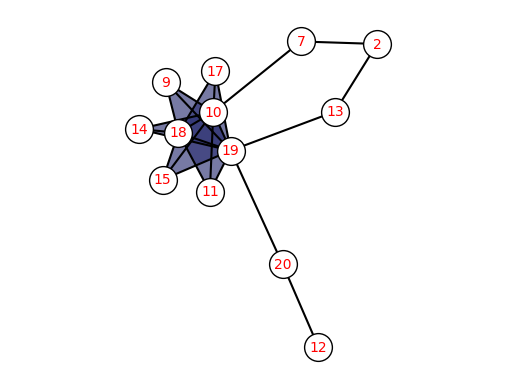}
		\caption{$\delta_{CC}=0.30,\beta_0=1$}
	\end{subfigure} \\
	\begin{subfigure}[b]{0.3\textwidth}
		\centering
		\includegraphics[width=\textwidth]{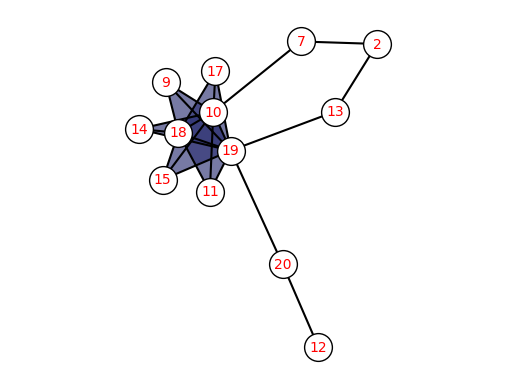}
		\caption{$\delta_{CC}=0.16,\beta_0=1$}
	\end{subfigure}
	\begin{subfigure}[b]{0.3\textwidth}
		\centering
		\includegraphics[width=\textwidth]{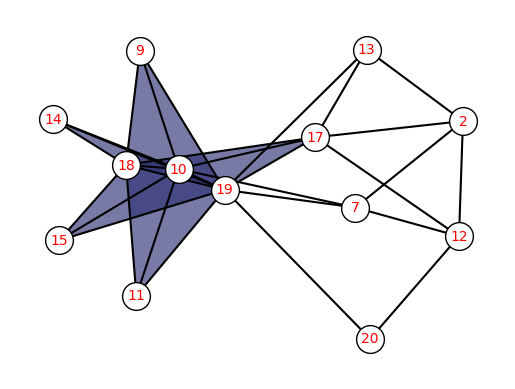}
		\caption{$\delta_{CC}=0.09,\beta_0=1$}
	\end{subfigure}
	\begin{subfigure}[b]{0.3\textwidth}
		\centering
		\includegraphics[width=\textwidth]{threshold_0.png}
		\caption{$\delta_{CC}=0, \beta_0=1$}
	\end{subfigure}
	\caption{Filtration of Lake network using the generalised clustering coefficient scores as threshold for $\delta_{CC}={1, 0.66, 0.30, 0.16, 0.09, 0}$}. 
	\label{fig:six_images}
\end{figure}

\section{Discussion}
In this work, we have used the idea of weighted networks to study the persistent homology on a given undirected, unweighted network. A novel adjacency rule based on the strength of adjacency has been given by which an unweighted network has converted into a weighted network, and its weighted simplicial adjacency matrices have been calculated at every structural level for its corresponding clique complex. These weighted simplicial adjacency matrices are then used to calculate the global measure, that is generalised betweenness centrality and using the scores obtained by this global measure a filtration scheme has been studied. This constructed filtration using the global measure reveals topological features of the studied network, which is the Lake network from the aquatic environment. Betti numbers reveals that there is absence of $1$ and $2$-dimensional holes in the studied clique complex for every threshold value except for $\delta_B$=0.02 and $0$.

One more generalised centrality, the maximal generalised degree, has been given in this work, which counts the degree of each simplex in a given simplicial complex. It tells us about the number of facets incident on a given simplex while considering their dimension as well. This is a local measure that is used for the filtration of simplicial complexes using degree scores. In both filtrations, we have observed a pattern that can be seen by the number of $1$-dimensional holes; it increases near the minimum value of centrality. As the set threshold reached near the minimum value, suddenly $1$ and $ 2$-dimensional holes appeared.

The variation in Betti numbers across both filtrations highlights that the topological features of a space shift similarly to how simplex rankings differ under various centrality indices. A simplex that holds significance from one perspective may lose relevance when the perspective changes. Likewise, topological features—such as holes—are influenced by adjustments in the filtration rule. Changes in Betti numbers for two distinct filtration schemes based on the two generalised centralities can be observed in Tables 2 and 3. In the Lake network, a notable pattern appears: the most significant shifts in Betti numbers occur near the minimum threshold value. For $\delta = 2$ and $\delta_{B_5} = 0.02$, a clear change in Betti numbers is noticeable across both filtrations.

Currently, the primary objective of most studies is to extract topological features of networks that persist across multiple scales. However, these studies face certain limitations, such that the computational difficulty associated with larger networks is a significant challenge. For instance, the network illustrated in this work is relatively small, consisting of only $20$ vertices.

\section*{Declaration of competing interest}
The authors declare that they have no known competing financial interests or personal relationships that could have appeared to influence the work reported in this paper.

\section*{Data availability}
We have included a link with all datasets used in the paper.

\section*{Funding} This research did not receive any specific grant from funding agencies in the public, commercial, or not-for-profit sectors.

\section*{Acknowledgements}
The authors, U.R and S.B., thank Shiv Nadar Institution of Eminence-Delhi NCR for providing facilities and resources for this research. The author S. M. acknowledges the support by the Ministry of Education, Science and Technological Development of the Republic of Serbia, with Contract Annex number 451-03-136/2025-03/200017, and by the Science Fund of the Republic of Serbia, 7416, Topology-derived methods for the analysis of collective trust dynamics - CTRUST.

Authors:\\

Udit Raj\\
Shiv Nadar Institution of Eminence, Delhi-NCR\\
ur376@snu.edu.in\\

Slobodan Maleti\'{c}\\
Vinča Institute of Nuclear Sciences-National Institute of the Republic of Serbia, University of Belgrade, Serbia\\
maletic.sloba@gmail.com\\

Sudeepto Bhattacharya\\
Shiv Nadar Institution of Eminence, Delhi-NCR\\
Sudeepto.bhattacharya@snu.edu.in


\begin{thebibliography}{}
		\bibitem{2} Aktas ME, Akbas E, Fatmaoui AE. Persistence homology of networks: methods and applications. Applied Network Science. 2019 Dec;4(1):1-28.
			
		\bibitem{Kartun-Giles2019}
		Alexander P. Kartun-Giles and Ginestra Bianconi, Beyond the clustering coefficient: A topological analysis of node neighbourhoods in complex networks, Chaos, Solitons \& Fractals: X, Volume 1, 2019, 100004
	\bibitem{9} Battiston F, Cencetti G, Iacopini I, Latora V, Lucas M, Patania A, Young JG, Petri G. Networks beyond pairwise interactions: Structure and dynamics. Physics reports. 2020 Aug 25;874:1-92.
	
	\bibitem{Boccaletti2006}
	Boccaletti S, Latora V, Moreno Y, Chavez M, Hwang DU. Complex networks: Structure and dynamics. Physics
	reports. 2006 Feb 1;424(4-5):175-308
	
	\bibitem{Boccaletti2023} Boccaletti S, De Lellis P, Del Genio CI, Alfaro-Bittner K, Criado R, Jalan S, Romance M. The structure and dynamics of networks with higher order interactions. Physics Reports. 2023 May 23;1018:1-64.
	
	\bibitem{16} Benson AR, Gleich DF, Leskovec J. Higher-order organization of complex networks. Science. 2016 Jul 8;353(6295):163-6.
	\bibitem{17} Bick C, Gross E, Harrington HA, Schaub MT. What are higher-order networks?. SIAM Review. 2023;65(3):686-731.
	\bibitem{1} Carstens CJ, Horadam KJ. Persistent homology of collaboration networks. Mathematical problems in engineering. 2013 Jun;2013.
	\bibitem{Bian}Courtney OT, Bianconi G. Generalized network structures: The configuration model and the canonical ensemble of simplicial complexes. Physical Review E. 2016 Jun 16;93(6):062311.
	\bibitem{7} Courtney O, Bianconi G. Weighted growing simplicial complexes. Phys Rev E 2017;95(6). http://dx.doi.org/10.1103/physreve.95.062301
		\bibitem{11} Deo S. Algebraic topology. Springer Singapore; 2018.
			\bibitem{Serrano2020}
		D. Hern\'{a}ndez Serrano and D. S\'{a}nchez G\'{o}mez, Centrality measures in simplicial complexes: Applications of topological data analysis to network science, Applied Mathematics and Computation 382 (2020) 125331
			\bibitem{13} Dunaeva O, Edelsbrunner H, Lukyanov A, Machin M, Malkova D, Kuvaev R, Kashin S. The classification of endoscopy images with persistent homology. Pattern Recognition Letters. 2016 Nov 1;83:13-22.
	\bibitem{15} Estrada E, Ross GJ. Centralities in simplicial complexes. Applications to protein interaction networks. Journal of theoretical biology. 2018 Feb 7;438:46-60.
	\bibitem{3} Hagberg A, Conway D. Networkx: Network analysis with python. URL: https://networkx. github. io. 2020.

	\bibitem{Munkres}
	Munkres JR. Elements of algebraic topology
	
	\bibitem{21} Matsumoto Y. An introduction to Morse theory. American Mathematical Soc.; 2002.
	\bibitem{19} Newman M. Networks. Oxford university press; 2018 Jul 4.
	\bibitem{12} Nanda V. Discrete Morse theory for filtrations (Doctoral dissertation, Rutgers University-Graduate School-New Brunswick).
	\bibitem{6} Otter N, Porter MA, Tillmann U, Grindrod P, Harrington HA. A roadmap for the computation of persistent homology. EPJ Data Science. 2017 Dec;6:1-38.
	
	\bibitem{5} Patania A, Vaccarino F, Petri G. Topological analysis of data. EPJ Data Science. 2017 Dec;6(1):1-6.
	\bibitem{Raj}Raj U, Banerjee A, Ray S, Bhattacharya S (2024) Structure of higher-order interactions in social-ecological networks through Q-analysis of their neighbourhood and clique complex. PLoS ONE 19(8): e0306409. https://doi.org/10.1371/journal.pone.0306409
	\bibitem{4} Raj U, Bhattacharya S. Some generalized centralities in higher-order networks represented by simplicial complexes. Journal of Complex Networks. 2023 Oct 1;11(5):cnad032.
	\bibitem{23} Robinson M. Hunting for oxes with sheaves 2019
	\bibitem{18} Strogatz S. Sync: The emerging science of spontaneous order.
    \bibitem{22} Saadatfar M, Takeuchi H, Robins V, Francois N, Hiraoka Y. Pore configuration landscape of granular crystallization. Nature communications. 2017 May 12;8(1):15082.
	\bibitem{14} Serrano DH, Hernández-Serrano J, Gómez DS. Simplicial degree in complex networks. Applications of topological data analysis to network science. Chaos, Solitons and Fractals. 2020 Aug 1;137:109839.
	\bibitem{Maletic2008}S. Maleti\' c , M. Rajkovi\' c , D. Vasiljevi\' c , Simplicial complexes of networks and their statistical properties,   in: M. Bubak, G.D. van Albada, J. Dongarra, P.M.A. Sloot (Eds.), Proceedings of the Computational Science ICCS, Lecture Notes in Computer Science, 5102, Springer, Berlin, Heidelberg, 2008	
\end{thebibliography}
\end{document}